\documentclass[preprint,12pt]{elsarticleb}

\usepackage{amssymb}
\usepackage{graphicx}
\usepackage{amsmath}
\usepackage{exscale}
\usepackage{times}
\usepackage{color}

\font\pppppcarac=ptmr8y at 5pt
\font\ppppcarac=ptmr8y at 6pt

\font\ppcarac=ptmr8y at 8pt

\font\carac=ptmr8y at 11pt

\font\bf=ptmb8y at 10pt

\newcommand{\bfH}{{\textbf{H}}}

\newcommand{\bfU}{{\textbf{U}}}
\newcommand{\bfV}{{\textbf{V}}}
\newcommand{\bfW}{{\textbf{W}}}
\newcommand{\bfX}{{\textbf{X}}}

\newcommand{\bfZ}{{\textbf{Z}}}

\newcommand{\bfg}{{\textbf{g}}}

\newcommand{\bfm}{{\textbf{m}}}

\newcommand{\bfu}{{\textbf{u}}}
\newcommand{\bfv}{{\textbf{v}}}

\newcommand{\bfx}{{\textbf{x}}}

\newcommand{\bfeta}{{\boldsymbol{\eta}}}

\newcommand{\bfvarphi}{{\boldsymbol{\varphi}}}
\newcommand{\bfpsi}{{\boldsymbol{\psi}}}

\newcommand{\bfphi}{{\boldsymbol{\phi}}}

\font\teneclair=bbold12
\font\seveneclair=bbold9
\font\fiveeclair=bbold7
\newfam\eclairfam
\scriptscriptfont\eclairfam=\fiveeclair
\textfont\eclairfam=\teneclair
\scriptfont\eclairfam=\seveneclair
\def\eclair{\fam\eclairfam\teneclair}

\def\11{{\eclair 1}}

\def\MM{{\eclair{M}}}

\def\PP{{\eclair{P}}}
\def\RR{{\eclair{R}}}

\newcommand{\curH}{{\mathcal{H}}}

\newcommand{\curL}{{\mathcal{L}}}

\newcommand{\curP}{{\mathcal{P}}}
\newcommand{\curR}{{\mathcal{R}}}
\newcommand{\curT}{{\mathcal{T}}}

\newcommand{\curV}{{\mathcal{V}}}

\newcommand{\bfcurL}{{\boldsymbol{\mathcal{L}}}}
\newcommand{\bfcurN}{{\boldsymbol{\mathcal{N}}}}
\newcommand{\bfcurW}{{\boldsymbol{\mathcal{W}}}}
\newcommand{\bfcurY}{{\boldsymbol{\mathcal{Y}}}}
\newcommand{\bfcurZ}{{\boldsymbol{\mathcal{Z}}}}

\newcommand{\bfzero}{{ \hbox{\bf 0} }}
\newcommand{\red}{{\hbox{{\ppppcarac red}}}}
\newcommand{\tr}{{\hbox{{\textrm tr}}}}
\newcommand{\Fac}{{\hbox{{\carac Fac}}}}
\newcommand{\curSn}{{{\mathcal{S}}_n}}
\newcommand{\curSnu}{{{\mathcal{S}}_\nu}}
\newcommand{\MC}{{\hbox{{\ppppcarac MC}}}}
\newcommand{\pMC}{{\hbox{{\pppppcarac MC}}}}
\newcommand{\st}{{\hbox{{\ppcarac st}}}}

\journal{Journal of Computational Physics on 27 September 2015}

\begin{document}

\begin{frontmatter}



\title{Data-driven probability concentration\\
 and sampling on manifold}

\author[1]{C. Soize \corref{cor1}}
\ead{christian.soize@univ-paris-est.fr}
\author[2]{R. Ghanem}
\ead{ghanem@usc.edu}
\cortext[cor1]{Corresponding author: C. Soize,
  christian.soize@univ-paris-est.fr}
\address[1]{Universit\'e Paris-Est, Laboratoire Mod\'elisation et
  Simulation Multi-Echelle, MSME UMR 8208 CNRS, 5 bd Descartes, 77454
  Marne-La-Vall\'ee, Cedex 2, France}
\address[2]{University of Southern California, 210 KAP Hall, Los
  Angeles, CA 90089, United States}

\date{27 September 2015}

\begin{abstract}
A new methodology is proposed for generating realizations of a random
vector with values in a finite-dimensional Euclidean space that are
statistically consistent with a data set of observations of this vector.
The probability distribution of this random vector, while a-priori not
known, is presumed to be concentrated on an unknown subset of the
Euclidean space. A random matrix is introduced whose columns are
independent copies of the random vector and for which the number of
columns is the number of data points in the data set. The approach is
based on the use of (i) the multidimensional kernel-density estimation
method for estimating the probability distribution of the random
matrix, (ii) a MCMC method for  generating realizations for the random
matrix, (iii) the diffusion-maps approach  for discovering and
characterizing the geometry and the structure of the data set, and
(iv) a reduced-order representation of the random matrix, which is
constructed using the diffusion-maps vectors associated with the first
eigenvalues of the transition matrix relative to the given data
set. The convergence aspects of the proposed methodology are
analyzed and a numerical validation is explored through three
applications of increasing complexity. The proposed method is found to
be robust to noise levels and data complexity as well as to the intrinsic
dimension of data and the size of experimental data sets.  Both the
methodology and the underlying mathematical framework presented in
this paper contribute new capabilities and perspectives at the
interface of uncertainty quantification, statistical data analysis,
stochastic modeling and associated statistical inverse problems.
\end{abstract}


\end{frontmatter}
\section{Introduction}
\label{Section1}
The construction of a generator of realizations from a given data set
related to a $\RR^n$-valued random vector, for which the probability
distribution is unknown and is concentrated on an unknown subset
$\curSn$ of $\RR^n$, is a central and difficult problem in uncertainty
quantification and statistical data analysis, in stochastic modeling
and associated statistical inverse problems for boundary value
problems, in the design of experiments for random parameters, and
certainly, in signal processing and machine learning.\\
Two fundamental tools serve as building blocks for addressing this problem.
First, nonparametric statistical methods \cite{Bowman1997,Scott2015}
can be effectively used to construct probability distribution on
$\RR^n$ of a random vector given an initial data set of its samples.
Multidimensional Gaussian kernel-density estimation is one efficient
subclass of these methods. Markov chain Monte Carlo (MCMC) procedures
can then be used to sample additional realizations from the resulting
probability model, and which are thus statistically consistent with
the initial data set \cite{Kaipio2005,Robert2005,Spall2003}.  The
second building block consists of manifold embedding algorithms, where
low-dimensional structure is characterized within a larger vector
space.   Diffusion maps \cite{Coifman2005,Coifman2006,Talmon2015} is a
powerful tool for characterizing and delineating $\curSn$ using the
initial data set and concepts of geometric diffusion.\\
The first tool described above, consisting of using nonparametric
density estimation with MCMC, does not allow, in general, the restriction
of new samples to the subset $\curSn$ on which the probability
distribution is concentrated.  The scatter of generated samples
outside of $\curSn$ is more pronounced the more complex and
disconnected this set is.  \\
The second tool consisting of diffusion maps, while effectively
allowing for the discovery and characterization of subset $\curSn$ on
which the probability distribution is concentrated, does not give a
direct approach for generating additional realizations in this subset
that are drawn from a target distribution consistent with the initial
data set.\\
These two fundamental tools have been used independently and quite
successfully to address problems of sampling from complex probability
models and detecting low-dimensional manifolds in high-dimensional
settings.  An analysis of MCMC methods on Riemann manifolds has been
presented recently \cite{Girolami2011} where the manifold is the locus
of density functions and not of the data itself.
This paper addresses the still open challenge of efficient statistical
sampling on manifolds defined by limited data.\\
It should be noted that the PCA \cite{Jolliffe2002} yields a statistical reduction method for second-order random vectors in finite dimension, similarly to the Karhunen-Lo\`eve expansion (KLE) \cite{Karhunen1946,Loeve1955}, which yields a statistical reduction method for second-order stochastic processes and random fields, and which has been used for obtaining an efficient construction \cite{Ghanem1990,Ghanem1991} of the polynomial chaos expansion (PCE) of stochastic processes and random fields \cite{Cameron1947}, and for which some ingredients have more recently been introduced for analyzing  complex problems encountered in uncertainty quantification \cite{Soize2004,Perrin2013}. \textit{A priori} and in general, the PCA or the KLE, which use a nonlocal basis with respect to the data set (global basis related to the covariance operator estimated with the data set) does not allow for discovering and characterizing the subset on which the probability law is concentrated.
The present work can be viewed as an extension and generalization of
previous work by the authors where the low-dimensional manifold
was unduly restricted \cite{Tipireddy2014,Ghanem2015,Thimmisetty2015}.\\
After formalizing the problem in Section~\ref{Section2}, the proposed
methodology is presented in Section~\ref{Section3} and developed in
Section~\ref{Section4}. Section~\ref{Section5} deals with three
applications: the first two applications correspond to analytical examples in dimension $2$ with $230$ given data points and in dimension $3$ with $400$ data points. The third application is related to a petro-physics database made up of experimental measurements for which the dimension is $35$ with $13,056$ given data points.
\subsection*{Notations}
\noindent A lower case letter such as $x$, $\eta$, or $u$, is a  real deterministic variable.\\
A boldface lower case letter such as $\bfx$, $\bfeta$, or $\bfu$ is a real deterministic vector.\\
An upper case letter such as $X$, $H$, or $U$, is a real random variable.\\
A boldface upper case letter, $\bfX$, $\bfH$, or $\bfU$, is a real random vector.\\
A lower case letter between brackets such as  $[x]$, $[\eta]$, or $[u]$), is a real deterministic matrix.\\
A boldface upper case letter between brackets such as $[\bfX]$, $[\bfH]$, or $[\bfU]$, is a real random matrix.\\

\noindent
$E$: Mathematical expectation.\\
$\MM_{n,N}$: set of all the $(n\times N)$ real matrices.\\
$\MM_\nu$: $\MM_{\nu,\nu}$.\\
$\Vert\bfx\Vert$: Euclidean norm of vector $\bfx$.\\
$[x]_{kj}$: entry of matrix $[x]$.\\
$[x]^T$: transpose of matrix $[x]$.\\
$\tr\{[x]\}$: trace of a square matrix $[x]$.\\
$\Vert [x ]\Vert_F$: Frobenius norm of matrix $[x]$ such that $\Vert x \Vert_F^2=\tr\{[x]^T\, [x]\}$.\\
$[I_{\nu}]$: identity matrix in $\MM_\nu$.\\
$\delta_{kk'}$: Kronecker's symbol such that $\delta_{kk'} =0$ if $k\not= k'$ and $=1$ if $k=k'$.\\
\section{Problem set-up}
\label{Section2}
The following four ingredients serve to set the stage for the
mathematical analysis required for constructing the target
probability distribution and sampling from it.

\noindent (i) Let $\bfx=(x_1,\ldots , x_n)$ be a generic point in $\RR^n$ and  let $d\bfx = dx_1\ldots dx_n$ be the Lebesgue measure.
 A family of $N$ vectors in $\RR^n$ will be written as $\{\bfx^1, \ldots , \bfx^N\}$.\\

\noindent (ii) Let $\bfX=(X_1,\ldots ,X_n)$ be a random vector defined
on a probability space $(\Theta,\curT,\curP)$, with values in $\RR^n$,
for which the probability distribution  is defined by a probability
density function (pdf) on $\RR^n$  (\textit{a priori} and in general,
the probability distribution is not Gaussian). This pdf is unknown but
is assumed to be concentrated on an unknown subset $\curSn$ of
$\RR^n$. A specific realization of random vector $\bfX$ will be
denoted by $\bfX(\theta)$ where $\theta\in\Theta$. \\

\noindent (iii) The available information consists of a given set of
$N$ data points specified by $N$ vectors $\bfx^{d,1},\ldots
,\bfx^{d,N}$ in $\RR^n$. These will be assumed to constitute $N$
statistically independent realizations (or samples)
$\bfX(\theta_1),$ $\ldots,\bfX(\theta_N)$ of random vector $\bfX$.
For $j=1,\ldots , N$, the vector $\bfx^{d,j}$ in $\RR^n$ is written as  $\bfx^{d,j} = (x^{d,j}_1,\ldots , x^{d,j}_n)$.
The $N$ data points can then be represented by the matrix $[x_d]$ in $\MM_{n,N}$ such that $[x_d]_{kj} = x^{d,j}_k$.\\

\noindent (iv) The local structure of the given data set is captured
via random matrix $[\bfX]$, defined on $(\Theta,\curT,\curP)$, with
values in $\MM_{n,N}$.  Specifically, $[\bfX]=[\bfX^1 \ldots \bfX^N]$ in
which the columns $\bfX^1, \ldots ,\bfX^N$ are $N$ independent copies
of random vector $\bfX$.  Consequently, matrix $[x_d]$ can be viewed
as one realization of random matrix $[\bfX]$\\

The objective of this paper then is to construct a generator
of realizations of random matrix $[\bfX]$ in $\MM_{n,N}$, for which
the unknown probability distribution is directly deduced from the
unknown probability distribution of random vector  $\bfX$, which  is
concentrated on the unknown subset $\curSn$ of $\RR^n$, and for which
only one realization $[x_d]$ is given.\\

The unknown subset $\curSn$ of $\RR^n$ can be viewed as a manifold,
which corresponds to the structure of data $[x_d]$,  and on which the
unknown probability measure is concentrated. Consequently, the
objective of the paper is to perform "data-driven probability concentration and sampling on a manifold".
\section{Summary of the methodology proposed}
\label{Section3}
To enhance the utility of the present paper and to clarify the
inter-relation between a number of intricate mathematical steps, the
proposed methodology is summarized in the following seven steps.
\begin{enumerate}
\item In general, the given data  are heterogeneous and badly conditioned. Consequently, the first step consists in performing a scaling of the given data, which yields the matrix $[x_d]$ in $\MM_{n,N}$ of the scaled given data (the matrix introduced in Section~\ref{Section2}), and simply called the given data set (removing the word "scaled"). The given data set are then normalized by using a principal component analysis (but without trying to introduce a statistical reduced-order representation). Therefore, the random matrix $[\bfX]$ (corresponding to scaled data $[x_d]$) is written as an affine transformation of a random matrix $[\bfH]$ with values in $\MM_{\nu,N}$ with $1 < \nu \leq n$ (in general, $\nu = n$, but sometimes some eigenvalues (of the empirical estimate of the covariance matrix of $\bfX$) exhibits zeros eigenvalues that are removed, yielding $\nu < n$). Random matrix $[\bfH]$ can then be written as $[\bfH] = [\bfH^1 \ldots \bfH^N]$ in which the columns $\bfH^1,\ldots ,\bfH^N$ are $N$ independent copies of a random vector $\bfH$ with values in $\RR^\nu$, whose probability density function on $\RR^\nu$ is unknown and is concentrated on an unknown subset $\curSnu$ of $\RR^\nu$. The given data $[x_d]$ in $\MM_{n,N}$ (related to $[\bfX]$) are then transformed into given data, $[\eta_d]$ in $\MM_{\nu,N}$, related to random matrix $[\bfH]$. The data represented by $[\eta_d]$ are thus normalized. Let $p_\bfH$ be the nonparametric estimate of the probability density function of random vector $\bfH$, which is performed by using $[\eta_d]$ (note that $p_\bfH$ is not the pdf of $\bfH$ but is the nonparametric estimate of the pdf of $\bfH$). Consequently,  the nonparametric estimate of the probability distribution on $\MM_{\nu,N}$ of random matrix $[\bfH]$ is written as $p_{[\bfH]}([\eta])\, d[\eta]= p_\bfH(\bfeta^1)\times \ldots \times p_\bfH (\bfeta^N) \, d\bfeta^1\ldots d\bfeta^N$ in which $[\eta]$ is any matrix in $\MM_{\nu,N}$ such that $[\eta] = [\bfeta^1 \ldots \bfeta^N]$ with $\bfeta^j\in\RR^\nu$.
\item The second step consists in constructing the nonparametric statistical estimate $p_\bfH$ of the probability density function of $\bfH$ using data $[\eta_d]\in\MM_{\nu,N}$. This is an usual problem that will be performed by using the classical multidimensional Gaussian kernel-density estimation method. Nevertheless, we will use the modification proposed in \cite{Soize2015} (instead of the classical method) in order that the nonparametric estimate $p_\bfH$ yields, for the estimation of the covariance matrix of $\bfH$ (using $[\eta_d]$), the identity matrix $[I_\nu]$ in $\MM_\nu$.
    This construction is directly used in the following third step.
\item The third step consists in introducing an adapted  generator of realizations for random matrix $[\bfH]$, which belongs to the class of the MCMC methods  such as the Metropolis-Hastings algorithm \cite{Metropolis1949,Hastings1970} (that requires the definition of a good proposal distribution), the Gibbs sampling \cite{Geman1984} (that requires the knowledge of the conditional distribution)  or the slice sampling \cite{Neal2003} (that can exhibit difficulties related to the general shape of the probability distribution, in particular for multimodal  distributions). This adapted generator will be the one derived from \cite{Soize2015}, which is based on a nonlinear It\^o stochastic differential equation (ISDE) formulated for a dissipative Hamiltonian dynamical system \cite{Soize2008}, which admits $p_{[\bfH]}([\eta])\, d[\eta]$ as an invariant measure, and for which the initial condition depends on matrix $[\eta_d]$.
\item The fourth step of the methodology consists in characterizing the subset $\curSnu$ from scaled and normalized data $[\eta_d]$. This will be done using the formulation of the diffusion maps, which is a very powerful mathematical tool for doing that. It should be noted that the diffusion-maps method is a local approach with respect to given data while the PCA is a global approach that, in general, cannot see the local geometric structure of the given data set. However, the  diffusion distance, which has been introduced in \cite{Coifman2005} for discovering and characterizing $\curSnu$, and which is constructed using the diffusion maps, does not allow for constructing a generator of realizations of random matrix $[\bfH]$ for which data $[\eta_d]$ are given but for which its probability measure and the subset $\curSnu$ of concentration are unknown. This step is introduced for constructing an algebraic vector basis $\{\bfg^1,\ldots,\bfg^N\}$ of $\RR^N$, depending on two parameters that are a smoothing parameter $\varepsilon > 0$ and an integer $\kappa$ related to the analysis scale of the local geometric structure of the data set. For $\alpha=1,\ldots , N$,  the vectors $\bfg^\alpha = (g^\alpha_1,\ldots,g^\alpha_N)\in \RR^N$ are directly related to the diffusion maps. A subset of this basis will be able to characterize the subset $\curSnu$ of $\RR^\nu$ on which the probability measure of $\bfH$ is concentrated. We will then introduce the matrix $[g]$ in $\MM_{N,m}$ made up of the first $m$ vectors $\{\bfg^1,\ldots,\bfg^m\}$ of the diffusion-maps basis, with $1 < m \ll N$.
\item  The fifth step  consists in estimating an adapted value of $m$ in order to capture the local geometric structure of $\curSnu$ and to obtain a reasonable mean-square convergence.
\item Using the first $m$ vectors (represented by matrix $[g]$) of the diffusion-maps basis, the sixth step consists in constructing a reduced-order ISDE, which allows for  generating some additional realizations of the reduced-order representation of random matrix $[\bfH]$, by introducing the random matrix $[\bfZ]$ with values in $\MM_{\nu,m}$ such that $[\bfH] = [\bfZ]\, [g]^T$.
\item The last step consists in numerically solving the reduced-order ISDE for computing the additional realizations $[z_s^1], \ldots , [z_s^{n_\pMC}]$ of random matrix $[\bfZ]$ and then to deduce the additional realizations $[x_s^1], \ldots , [x_s^{n_\pMC}]$ of random matrix $[\bfX]$ for which only one realization $[x_d]$ was given.
\end{enumerate}
\section{Formulation}
\label{Section4}
In this section, a detailed presentation of the methodology is given that parallels
the steps described in Section~\ref{Section3}.
\subsection{Scaling and normalizing the given data set}
\label{Section4.1}
Let $[x_d^{uns}]$ be the matrix in $\MM_{n,N}$ of the unscaled given data set. The matrix  $[x_d]$ in $\MM_{n,N}$ of the scaled given data set
(simply called the given data set) is constructed (if the data effectively require such a scaling, which will be the case for the third
application presented in Section~\ref{Section5.3}) such that, for all $k=1,\ldots , n$ and $j=1,\ldots , N$,
\begin{equation}
\label{EQ1}
[x_d]_{kj} = \frac{[x_d^{uns}]_{kj} - \min_{j'} [x_d^{uns}]_{kj'} }{\max_{j'} [x_d^{uns}]_{kj'} - \min_{j'} [x_d^{uns}]_{kj'} } + \epsilon_s\, .   
\end{equation}
The quantity $\epsilon_s$ is added to the scaled data in order to avoid the scalar $0$ in the nonparametric statistical estimation of the pdf.
Let $\bfm$ and $[c]$ be the empirical estimates of the mean vector $E\{\bfX\}$ and the covariance matrix
$E\{(\bfX -E\{\bfX\})\, (\bfX -E\{\bfX\})^T\}$, such that
\begin{equation}
\label{EQ2}
\bfm = \frac{1}{N} \sum_{j=1}^N \bfx^{d,j}\quad , \quad
 [c] = \frac{1}{N-1} \sum_{j=1}^N (\bfx^{d,j}-\bfm)\, (\bfx^{d,j} -\bfm)^T \, .                                                                  
\end{equation}
We consider the eigenvalue problem $[c]\, \bfvarphi^k = \mu_k\,
\bfvarphi^k$. Noting that matrix $[c]$ is often of rank $\nu \leq n$,
denote its $\nu$ positive eigenvalues by $\{\mu_i\}_{i=1}^\nu$ with $0 < \mu_1\leq \mu_2\leq \ldots \leq \mu_\nu$ and let $[\varphi]$ be the $(n\times\nu)$  matrix such $[\varphi]^T\,[\varphi]= [I_\nu]$, whose columns are the associated orthonormal eigenvectors $\bfvarphi^1,\ldots , \bfvarphi^\nu$.
 Consequently, random matrix $[\bfX]$ can be rewritten as
\begin{equation}
\label{EQ3}
[\bfX] = [\underline x] + [\varphi]\, [\mu]^{1/2}\, [\bfH] \, ,                                                                                   
\end{equation}
in which $[\underline x]$ is the matrix in $\MM_{n,N}$ for which each column is vector $\bfm$ and  where $[\mu]$ is the positive diagonal $(\nu\times\nu)$ real matrix such that $[\mu]_{kk'}  = \delta_{kk'}\mu_k$. The realization $[\eta_d]\in \MM_{\nu,N}$ of $[\bfH]$ associated with the realization $[x_d]$ of $[\bfX]$ is thus computed by
\begin{equation}
\label{EQ4}
[\eta_d] =  [\mu]^{-1/2} [\varphi]^T\, ([x_d] - [\underline x]) \, .                                                                              
\end{equation}
Let $\bfeta^{d,1},\ldots ,\bfeta^{d,N}$ be the $N$ vectors in $\RR^\nu$ such that
$[\bfeta^{d,1} \ldots \bfeta^{d,N}] = [\eta_d]$ (the columns of $[\eta_d]$).
It can easily be seen that the  empirical estimates $\bfm'$ of the mean vector $E\{\bfH\}$ and $[c']$ of the covariance matrix
$E\{(\bfH -E\{\bfH\})\, (\bfH -E\{\bfH\})^T\}$ of random vector $\bfH$ are such that
\begin{equation}
\label{EQ5}
\bfm' = \frac{1}{N} \sum_{j=1}^N \bfeta^{d,j} = \bfzero
\quad , \quad   [c'] = \frac{1}{N-1} \sum_{j=1}^N \bfeta^{d,j} \, (\bfeta^{d,j})^T = [\, I_\nu] \, .                                             
\end{equation}
\subsection{Construction of a nonparametric estimate $p_\bfH$ of the pdf of $\bfH$}
\label{Section4.2}
The estimation $p_\bfH$ on $\RR^\nu$ of the pdf of random vector $\bfH$ is carried out by using the Gaussian kernel-density estimation method and the $N$ independent realizations  $\bfeta^{d,1},\ldots , \bfeta^{d,N}$ represented by matrix $[\eta_d]$ computed with Eq.~\eqref{EQ4}.
As proposed in \cite{Soize2015}, a modification of the classical Gaussian kernel-density estimation method  is used  in order that the mean vector and the covariance matrix (computed with the nonparametric estimate $p_{{}_\bfH}$) are equal to $\bfzero$ and $[\, I_\nu]$ respectively (see Eq.~\eqref{EQ5}). The positive-valued function $p_\bfH$ on $\RR^\nu$ is then defined, for all $\bfeta$ in $\RR^\nu$,  by
\begin{equation}
\label{EQ6}
 p_\bfH(\bfeta) = \frac{1}{N}\sum_{j=1}^N \pi_{\nu,\widehat s_\nu}\,( \frac{\widehat s_\nu}{s_\nu} \bfeta^{d,j}-\bfeta )  \, ,            
\end{equation}
in which $\pi_{\nu,\widehat s_\nu}$ is the positive function from $\RR^\nu$ into $]0\, , +\infty[$ defined, for all $\bfeta$ in $\RR^\nu$,  by
\begin{equation}
\label{EQ7}
 \pi_{\nu,\widehat s_\nu}(\bfeta) = \frac{1}{(\sqrt{2\pi} \, \widehat s_\nu\,)^\nu}
           \, \exp\{- \frac{1}{2 {\widehat s_\nu}^{\, 2}} \Vert\bfeta\Vert^2\}   \, ,                                                     
\end{equation}
with $\Vert\bfeta\Vert^2 = \eta_1^2 +\ldots +\eta_\nu^2$ and where the positive parameters $s_\nu$ and $\widehat s_\nu$ are defined by
\begin{equation}
\label{EQ8}
s_\nu = \left\{\frac{4}{N(2+\nu)} \right\}^{1/(\nu+4)}   \quad , \quad  \widehat s_\nu =   \frac{s_\nu}{\sqrt{s_\nu^2 +\frac{N -1}{N}}} \, .   
\end{equation}
Parameter $s_\nu$ is the usual multidimensional optimal Silverman bandwidth (in taking into account that the empirical estimate of the standard deviation of each component is unity), and parameter $\widehat s_\nu$ has been introduced in order that the second equation in Eq.~\eqref{EQ5} holds.
Using Eqs.~\eqref{EQ6} to \eqref{EQ8}, it can easily be verified that
\begin{equation}
\label{EQ9}
\int_{\RR^\nu} \bfeta\, p_\bfH(\bfeta)\, d\bfeta = \frac{\widehat s_\nu}{s_\nu}\, \bfm'  =\bfzero \, ,              
\end{equation}
\begin{equation}
\label{EQ10}
\int_{\RR^\nu} \bfeta\,\bfeta^T \, p_\bfH(\bfeta)\, d\bfeta =
{\widehat s_\nu}^{\, 2}\, [\, I_\nu] + (\frac{\widehat s_\nu}{s_\nu})^2 \, \frac{(N - 1)}{N}[c'] = [\, I_\nu] \, .                  
\end{equation}
A nonparametric estimate $p_{[\bfH]}$ on $\MM_{\nu,N}$ of the probability density function of random matrix $[\bfH]$ is then written as
\begin{equation}
\label{EQ11}
p_{[\bfH]}([\eta])=p_\bfH(\bfeta^{d,1})\times \ldots\times p_\bfH(\bfeta^{d,N})\, ,                                                         
\end{equation}
in which $p_\bfH$ is defined by Eqs.~\eqref{EQ6} to \eqref{EQ8}.
\subsection{Construction of an ISDE for generating realizations of random matrix $[\bfH]$}
\label{Section4.3}
The probability density function defined by Eqs.~\eqref{EQ6} to \eqref{EQ8} is directly used for constructing the It\^o stochastic differential equation.
Let $\{ ([\bfU(r)],[\bfV(r)]), r\in \RR^+ \}$ be the Markov stochastic process defined on the probability space $(\Theta,\curT,$ $\curP)$, indexed
by $\RR^+= [0\, ,+\infty[$, with values in $\MM_{\nu,N}\times\MM_{\nu,N}$, satisfying, for all $r>0$,  the following ISDE
\begin{equation}
\label{EQ12}
 d[\bfU(r)] =  [\bfV(r)] \, dr \, ,                                                                                                          
\end{equation}
\begin{equation}
\label{EQ13}
d[\bfV(r)]=  [L([\bfU(r)])]\, dr -\frac{1}{2} f_0\, [\bfV(r)]\, dr + \sqrt{f_0}\, [d\bfW(r)] \, ,                                            
\end{equation}
with the initial condition
\begin{equation}
\label{EQ14}
[\bfU(0)] = [\bfH_d] \quad , \quad [\bfV(0)] = [\bfcurN\,] \quad a.s \, .                                                                    
\end{equation}
In Eqs.~\eqref{EQ13} and \eqref{EQ14}, the different quantities are defined as follows.\\

\noindent (i) For all $[u] = [\bfu^1 \ldots \bfu^N]$ in $\MM_{\nu,N}$ with $\bfu^\ell=(u^\ell_1,\ldots ,u^\ell_\nu)$ in $\RR^\nu$, the matrix $[L([u])]$ in $\MM_{\nu,N}$ is defined, for all $k = 1,\ldots ,\nu$ and for all $\ell=1,\ldots , N$, by
\begin{equation}
\label{EQ15}
[L([u])]_{k\ell}= -\frac{\partial}{\partial u^\ell_k}  \curV(\bfu^\ell) \, ,                                                                   
\end{equation}
in which the potential $\curV(\bfu^\ell)$ defined on $\RR^\nu$ with values in $\RR^+$, is defined by
\begin{equation}
\label{EQ16}
\curV(\bfu^\ell) = -\log\{q(\bfu^\ell)\} \, ,                                                                                                  
\end{equation}
where $\bfu^\ell\mapsto q(\bfu^\ell)$ is the continuously differentiable function from $\RR^\nu$ into $]0\, ,+\infty[$ such that
\begin{equation}
\label{EQ17}
q(\bfu^\ell) = \frac{1}{N} \sum_{j=1}^N
\exp\{ -\frac{1}{2 {\widehat s_\nu}^{\, 2}}\Vert \frac{\widehat s_\nu}{s_\nu}\bfeta^{d,j}-\bfu^\ell\Vert^2 \}  \, .                           
\end{equation}
From Eqs.~\eqref{EQ16} and \eqref{EQ17}, it can be deduced that,
\begin{equation}
\label{EQ17bis}
[L([u])]_{k\ell} = \frac{1}{q(\bfu^\ell)} \, \{ {\boldsymbol{\nabla}}_{\!\!\bfu^\ell}\, q(\bfu^\ell) \}_k \, ,                              
\end{equation}
\begin{equation}
\label{EQ17ter}
{\boldsymbol{\nabla}}_{\!\!\bfu^\ell}\, q(\bfu^\ell) = \frac{1}{\widehat s_\nu^{\,2}}\frac{1}{N} \sum_{j=1}^N  (\frac{\widehat s_\nu}{s_\nu}\bfeta^{d,j}-\bfu^\ell)\,
             \exp\{ -\frac{1}{2 \widehat s_\nu^{\,2}}\Vert \frac{\widehat s_\nu}{s_\nu}\bfeta^{d,j}-\bfu^\ell\Vert^2 \}  \, .                
\end{equation}
\\

\noindent (ii) The stochastic process $\{[d\bfW(r)], r \geq 0\}$ with values in $\MM_{\nu,N}$ is such that
$[d\bfW(r)] = [d\bfW^1(r) \ldots d\bfW^N(r)]$ in which the columns $\bfW^1 \ldots \bfW^N$ are $N$ independent copies of the normalized Wiener process $\bfW=(W_1,\ldots ,W_\nu)$ defined on $(\Theta,\curT,\curP)$, indexed by $\RR^+$ with values in $\RR^\nu$. The matrix-valued autocorrelation  function $[R_\bfW(r,r')]= E\{\bfW(r)\,\bfW(r')^T\}$ of $\bfW$ is then written as $[R_\bfW(r,r')] = \min (r,r')\, [I_\nu ]$.\\

\noindent (iii) The probability distribution of the random matrix $[\bfH_d]$ with values in $\MM_{\nu,N}$ is $p_{[\bfH]}([\eta])\, d[\eta]$. A known realization of $[\bfH_d]$ is matrix $[\eta_d]$. The random matrix $[\bfcurN\,]$ with values in $\MM_{\nu,N}$ is written as $[\bfcurN\,] = [\bfcurN^1 \ldots \bfcurN^N]$ in which the columns $\bfcurN^1, \ldots ,\bfcurN^N$ are $N$ independent copies of the normalized Gaussian vector $\bfcurN$ with values in $\RR^\nu$ (this means that $E\{\bfcurN\} =\bfzero$ and $E\{\bfcurN \bfcurN^T\} =[I_\nu]$). The random matrices $[\bfH_d]$ and $[\bfcurN\,]$, and the normalized Wiener
process $\{\bfW(r),r \geq 0\}$  are assumed to be independent.\\

\noindent (iv) The free parameter $f_0 > 0$ allows the dissipation term of the nonlinear second-order dynamical system (dissipative Hamiltonian system)  to be controlled.\\

Since the columns $\bfH^1,\ldots , \bfH^N$ of random matrix $[\bfH]$ are independent copies of random vector $\bfH$, and since the pdf of random matrix
$[\bfH_d]$ is $p_{[\bfH]}$, using Theorems 4 to 7 in pages 211 to 216 of Ref. \cite{Soize1994},
in which the Hamiltonian is taken as
$\curH(\bfu,\bfv) = \Vert \bfv\Vert^2 /2 + \curV(\bfu)$, and using \cite{Doob1990,Khasminskii2012} for proving the ergodic property, it can be proved  that the problem defined by Eqs.~\eqref{EQ12} to  \eqref{EQ14} admits a unique invariant measure and a unique solution $\{ ([\bfU(r)],[\bfV(r)]),$ $r\in \RR^+ \}$ that is a second-order diffusion stochastic process, which is stationary (for the shift semi-group on $\RR^+$ defined by the
positive shifts $r\mapsto r+ \tau$, $\tau \geq 0$) and ergodic, and such that, for all $r$ fixed in $\RR^+$,  the probability distribution of random matrix $[\bfU(r)]$ is $p_{[\bfH]}([\eta])\, d[\eta]$ in which  $p_{[\bfH]}$ is defined by Eq.~\eqref{EQ11}.\\

\noindent \textbf{Remarks}.\\

\noindent 1. It should be noted that the invariant measure is independent of $f_0$.\\

\noindent 2. If the initial condition $[\bfU(0)]$ was not $[\bfH_d]$ but was any other random matrix whose pdf is not $p_{[\bfH]}$, then the unique diffusion process $\{ ([\bfU(r)],[\bfV(r)]),$ $r\in \RR^+$ would not be stationary, but would be asymptotic (for $r\rightarrow +\infty$) to a stationary diffusion process $\{ ([\bfU_\st(r_\st)],$ $[\bfV_\st(r_\st)]), r_\st\geq 0\}$ such that, for all $r_\st > 0$,
$[\bfH] = [\bfU_\st(r_\st)] = \lim_{r\rightarrow +\infty} [\bfU(r)]$ in probability distribution (this implies that, for all $r_\st > 0$, the pdf of random matrix $[\bfU_\st(r_\st)]$ is $p_{[\bfH]}$). In such a case, the free parameter $f_0 > 0$ allows the transient response generated by the initial condition to be rapidly killed in order to get more rapidly the asymptotic behavior corresponding to the stationary and ergodic solution associated with the invariant measure.\\

\noindent 3. As the nonparametric estimate $p_{[\bfH]}$ of the pdf of $[\bfH]$ does not explicitly take into account the local structure of data set $[\eta_d]$, if the pdf of $\bfH$ is concentrated on $\curSnu$, then the generator of realizations constructed by the MCMC method defined by  Eqs.~\eqref{EQ12} to  \eqref{EQ14} (or by any other MCMC method), will not give some realizations localized in the subset $\curSnu$ (see the applications in Section~\ref{Section5}).\\

\noindent 4. As explained in \cite{Soize2015}, a variant of Eq.~\eqref{EQ13} could be introduced in replacing it by
$d[\bfV(r)]=  [L([\bfU(r)])]\, dr -\frac{1}{2} f_0\, [D_0]\,[\bfV(r)]\, dr + \sqrt{f_0}\,[S_0]\, [d\bfW(r)]$
in which $[S_0]$ would belong to $\MM_\nu$ and where $[D_0]$ would be a positive symmetric matrix such that $[D_0]=[S_0]\,[S_0]^T$ with
$1\leq \hbox{rank}[D_0] \leq \nu$. In the present case, such an extension would not allow for improving the methodology proposed because the initial condition for $[\bfU(0)]$ is the given matrix $[\eta_d]$ that follows $p_{[\bfH]}$.\\

\noindent 5. For $\theta$ fixed in $\Theta$, let $\{[\bfW(r;\theta)],r\geq 0\}$, $[\bfH_d(\theta)] = [\eta_d]$, and $[\bfcurN(\theta)]$  be independent realizations of the stochastic process $\{[\bfW(r)],r\geq 0\}$, the random matrix $[\bfH_d]$, and the random matrix $[\bfcurN]$.
Let $\{ ([\bfU(r;\theta)],[\bfV(r;\theta)]), r\in \RR^+\}$ be the corresponding realization of the unique stationary diffusion process $\{ ([\bfU(r)],[\bfV(r)]), r\in \RR^+\}$ of the ISDE problem defined by Eqs.~\eqref{EQ12} to  \eqref{EQ14}). Then additional realizations $[\eta_s^1], \ldots , [\eta_s^{n_\pMC}]$ of random matrix $[\bfH]$ can be generated by
\begin{equation}
\label{EQ18}
[\eta_s^\ell]  = [\bfU(\ell \rho ;\theta)] \quad , \quad \rho = M_0\, \Delta r \quad , \quad \ell = 1,\ldots , n_\pMC\, ,               
\end{equation}
in which $\Delta r$ is the sampling step of the continuous index parameter $r$ used in the integration scheme  (see Section~\ref{Section4.7.1})  and where $M_0$ is a positive integer:
\begin{itemize}
\item If $M_0 = 1$, then $\rho = \Delta r$ and the $n_\pMC$ additional realizations are dependent, but the ergodic property of $\{ [\bfU(r)],r\in \RR^+\}$ can be used for obtaining the convergence of statistics  constructed  using $[\eta_s^1], \ldots , [\eta_s^{n_\pMC}]$ for random matrix $[\bfH]$.
\item If integer $M_0$ is chosen sufficiently large (such that $\rho$
  is much larger than the relaxation time of the dissipative
  Hamiltonian dynamical system), then $[\eta_s^1], \ldots ,
  [\eta_s^{n_\pMC}]$ can approximatively be considered as independent
  realizations of random matrix $[\bfH]$.  We underscore here that
  each sample of random matrix $[\eta_s]$ consists of $N$ simultaneous
  samples of random vector $\bfH$ inherit additional statistical properties
  from the matrix structure of $[\bfH]$ to ensure their coalescence
  around the low-dimensional structure $\curSn$.
\end{itemize}
\subsection{Construction of a diffusion-maps basis $[g]$}
\label{Section4.4}
Let $k_\varepsilon(\bfeta,\bfeta')$ be the kernel defined on  $\RR^\nu\times \RR^\nu$, depending on a real smoothing parameter $\varepsilon > 0$, which verifies the following properties:
\begin{itemize}
\item $k_\varepsilon(\bfeta,\bfeta') = k_\varepsilon(\bfeta',\bfeta)$ (symmetry).
\item $k_\varepsilon(\bfeta,\bfeta')  \geq 0$ (positivity preserving).
\item $k_\varepsilon$ is positive semi-definite.
\end{itemize}
A classical choice (that we will use in Section~\ref{Section5}) for
the kernel that satisfies the above three properties is the Gaussian
kernel specified as,
\begin{equation}
\label{EQ19}
k_\varepsilon(\bfeta,\bfeta')  = \exp(-\frac{1}{4\varepsilon} \Vert\bfeta -\bfeta'\Vert^2) \, .                                              
\end{equation}
Let $[K]$ be the symmetric matrix in $\MM_N$ with positive entries such that
\begin{equation}
\label{EQ20}
[K]_{ij} = k_\varepsilon(\bfeta^{d,i}, \bfeta^{d,j}) \quad , \quad i \,\,\hbox{and} \,\, j \in \{1,\ldots , N\}\, .                        
\end{equation}
Let $[b]$ be the positive-definite diagonal real matrix in $\MM_N$ such that
\begin{equation}
\label{EQ21}
[b]_{ij} = \delta_{ij}\,\sum_{j'=1}^N [K]_{ij'} \, ,                                                                                        
\end{equation}
and let $[\PP]$ be the matrix in $\MM_N$ such that
\begin{equation}
\label{EQ22}
[\PP] = [b]^{-1}\, [K] \, .                                                                                                                 
\end{equation}
Consequently, matrix $[\PP]$ has positive entries and satisfies
$\sum_{j=1}^N [\PP]_{ij} = 1$ for all $i=1,\ldots , N$.  It can thus
be viewed as the transition matrix of a Markov chain that yields the probability of transition in one step. Let $[\PP_S]$ be the symmetric matrix in $\MM_N$ such that
\begin{equation}
\label{EQ23}
[\PP_S] = [b]^{1/2}\, [\PP] \, [b]^{-1/2} = [b]^{-1/2}\, [K] \, [b]^{-1/2} \, .                                                              
\end{equation}
We consider the eigenvalue problem $[\PP_S]\, \bfphi^\alpha = \lambda_\alpha\, \bfphi^\alpha$. Let $m$ be an integer such that $1 < m \leq N$.
It can easily be proved that the associated  eigenvalues are real, positive, and such that
\begin{equation}
\label{EQ24}
1= \lambda_1 > \lambda_2 \geq  \ldots \geq \lambda_m \, .                                                                  
\end{equation}
Let $[\phi]$ be the matrix in $\MM_{N,m}$ such that $[\phi]^T\, [\phi] = [I_m]$, whose columns are the $m$ orthonormal eigenvectors
$\bfphi^1,\ldots , \bfphi^m$ associated with $\lambda_1,\ldots , \lambda_m$. The eigenvalues of matrix $[\PP]$
are the same as the eigenvalues of matrix $[\PP_S]$. The right eigenvectors $\bfpsi^1,\ldots ,\bfpsi^m$ associated with
$\lambda_1,\ldots , \lambda_m$, which are such that $[\PP]\, \bfpsi^\alpha = \lambda_\alpha\, \bfpsi^\alpha$, are written as
\begin{equation}
\label{EQ25}
\bfpsi^\alpha = [b]^{-1/2} \, \bfphi^\alpha \in \RR^N \quad , \quad \alpha = 1,\ldots , m \, ,                                              
\end{equation}
and consequently, the matrix $[\psi] = [\bfpsi^1 \ldots \bfpsi^m] = [b]^{-1/2} \, [\phi] \in\MM_{N,m}$ is such that
\begin{equation}
\label{EQ26}
[\psi]^T\, [b]\, [\psi] = [I_m] \, ,                                                                                                          
\end{equation}
which defines the normalization of the right eigenvectors of $[\PP]$.\\

We then define a "diffusion-maps basis" by $[g] = [\bfg^{1} \ldots \bfg^m]\in \MM_{N,m}$ (which is an algebraic basis of $\RR^N$ for $m=N$) such that
\begin{equation}
\label{EQ27}
\bfg^\alpha =  \lambda_\alpha^\kappa\, \bfpsi^\alpha  \in \RR^N \quad , \quad \alpha=1,\ldots , m  \, ,                             
\end{equation}
in which  $\kappa$ is an integer that is chosen for fixing the analysis scale of the local geometric structure of the data set.
It should be noted that the family $\{\Psi_\kappa\}_\kappa$ of diffusion maps are defined \cite{Coifman2005,Coifman2006} by the vector $\Psi_\kappa =(\lambda_1^\kappa\, \bfpsi^1, \ldots, \lambda_m^\kappa\, \bfpsi^m)$ in order to construct a diffusion distance, and integer $\kappa$ is thus such that the probability of transition is in $\kappa$ steps. However, as we have previously explained, we do not use such a diffusion distance, but we use the "diffusion-maps basis"  $\{\bfg^{1} \ldots \bfg^N\}$ that we have introduced for performing a projection of each column of the $\MM_{N,\nu}$-valued random matrix $[\bfH]^T$ on the subspace of $\RR^N$, spanned by $\{\bfg^{1} \ldots \bfg^m\}$. Introducing the random matrix $[\bfZ]$ with values in $\MM_{\nu,m}$, we can then construct the following reduced-order representation of $[\bfH]$,
\begin{equation}
\label{EQ28}
[\bfH] = [\bfZ]\, [g]^T \, .                                                                                                             
\end{equation}
Since the matrix $[g]^T\, [g] \in \MM_m$ is invertible, Eq.~\eqref{EQ28} yields
\begin{equation}
\label{EQ29}
[\bfZ] = [\bfH]\, [a] \quad , \quad [a]  = [g]\, ([g]^T\, [g])^{-1} \in \MM_{N,m}\, .                                                      
\end{equation}
In particular, matrix $[\eta_d]\in \MM_{\nu,N}$ can be written as $[\eta_d] = [z_d]\, [g]^T$ in which the matrix $[z_d]\in \MM_{\nu,m}$ is written as
\begin{equation}
\label{EQ30}
[z_d] = [\eta_d]\, [a]  \in \MM_{\nu,m}\, .                                                                                                
\end{equation}
\subsection{Estimating dimension $m$ of the reduced-order representation of random matrix $[\bfH]$}
\label{Section4.5}
Because an estimation of the value of the order-reduction dimension $m$ must be known before beginning the generation of additional realizations of random matrix $[\bfZ]$ using the reduced-order representation of random matrix $[\bfH]$, we propose a methodology which is only based on the use of the known data set represented by matrix $[\eta_d]$ that is a realization of random matrix $[\bfH]$.\\

For a given value of integer $\kappa$ and for a given value of
smoothing parameter $\varepsilon > 0$, the decay of the graph
$\alpha\mapsto \lambda_\alpha$ of  the eigenvalues  of transition matrix $[\PP]$, yields a criterion for choosing
the value of $m$ that allows the local geometric structure of the data
set represented by  $[\eta_d]$ to be discovered. Nevertheless, this
criterion can be misleading as it does not capture statistical
fluctuations around the embedded manifold. An additional mean-square
convergence  must be verified, and if necessary, the value of $m$ must
be increased. However, if the value of $m$ is chosen too large, the
localization of the geometric structure of the data set is lost. Consequently, a compromise must be applied between the very small value of $m$ given by the decreasing criteria of the eigenvalues of matrix $[\PP]\in\MM_N$ and a larger value of $m$ which is necessary for obtaining a reasonable mean-square convergence.\\

Using Eqs.~\eqref{EQ28} to \eqref{EQ30} allows for calculating  the reduced-order representation $[\eta_\red(m)]\in\MM_{\nu,N}$ of $[\eta_d]$ such that
$[\eta_\red(m)] =[\eta_d]\,[a]\,[g]^T$ in which $[a]$ and $[g]$ depend on $m$. It should be noted that if $m=N$, then $[a]\,[g]^T = [I_N]$ and therefore,
$[\eta_\red(m)] =[\eta_d]$. In such a case, the "reduced-order" representation would correspond to a simple change of vector basis in $\RR^N$ and the localization of the geometric structure of the data set would be lost. This implies that $m$ must be much more less than $N$ for preserving the capability of the approach to localize the geometric structure of the data set, and must be chosen as the smallest possible value that yields a reasonable mean-square convergence.
Let $[x_\red(m)]\in\MM_{n,N}$ be the matrix $[x_d]$ of the data set, calculated using Eq.~\eqref{EQ3} with $[\eta_\red(m)]$. We then have
\begin{equation}
\label{EQ31}
[x_\red(m)] = [\underline x] + [\varphi]\, [\mu]^{1/2}\, [\eta_d]\,[a]\,[g]^T \, .                                                              
\end{equation}
Let $\bfx_\red^1(m),\ldots,\bfx_\red^N(m)$ be the $N$ vectors in $\RR^n$, which constitute the columns of matrix $[x_\red(m)]\in\MM_{n,N}$. We then introduced the empirical estimates $\bfm_\red\in\RR^N$ and $[c_\red]\in\MM_N$ of the mean value and the covariance matrix calculated with the realization
$[x_\red(m)]\in\MM_{n,N}$ such that
\begin{equation}
\label{EQ32}
\bfm_\red(m) = \frac{1}{N} \sum_{j=1}^N \bfx_\red^j(m) \, ,                                                                                      
\end{equation}
\begin{equation}
\label{EQ33}
[c_\red(m)] = \frac{1}{N-1} \sum_{j=1}^N (\bfx_\red^j(m)-\bfm_\red)\, (\bfx_\red^j(m) -\bfm_\red)^T \, .                                         
\end{equation}
The mean-square convergence criterion is then defined by
\begin{equation}
\label{EQ34}
e_\red(m) = \frac{\Vert [c_\red(m)] - [c] \Vert_F}{\Vert [c] \Vert_F} \,  .                                                            
\end{equation}
in which $[c]$ is defined by Eq.~\eqref{EQ2}. Since $[x_\red(N)] = [x_d]$, it can be deduced that $e_\red(m) \rightarrow 0$ when $m$ goes to $N$.
For a fixed reasonable value  $\epsilon_0 > 0$ of the relative tolerance $e_\red(m)$, an estimate of $m$ will consist in looking for the smallest value of $m$ such that $e_\red(m) \leq \varepsilon_0$.
An illustration of the use of this criterion will be given in the third application presented in Section~\ref{Section5.3}.
\subsection{Reduced-order ISDE for generation of additional realizations of random matrix $[\bfX]$}
\label{Section4.6}
For $m$, $\varepsilon$, and $\kappa$ fixed, the reduced-order representation  $[\bfH] = [\bfZ]\, [g]^T$ of random matrix $[\bfH]$, defined by Eq.~\eqref{EQ28}, is used for constructing the reduced-order ISDE associated with Eqs.~\eqref{EQ12} to \eqref{EQ14}.  Introducing the change of stochastic processes $[\bfU(r)]=[\bfcurZ(r)]\,[g]^T$ and $[\bfV(r)] = [\bfcurY(r)]\,[g]^T$ into  these equations, then right multiplying the obtained equations by matrix $[a]$, and taking into account Eq.~\eqref{EQ29}, it can be seen that $\{ ([\bfcurZ(r)],[\bfcurY(r)]), r\in \RR^+ \}$ is a Markov stochastic process defined on the probability space $(\Theta,\curT,$ $\curP)$, indexed by $\RR^+= [0\, ,+\infty[$, with values in $\MM_{\nu,m}\times\MM_{\nu,m}$, satisfying, for all $r>0$,  the following reduced-order ISDE,
\begin{equation}
\label{EQ35}
 d[\bfcurZ(r)] =  [\bfcurY(r)] \, dr \, ,                                                                                                     
\end{equation}
\begin{equation}
\label{EQ36}
d[\bfcurY(r)]=  [\curL([\bfcurZ(r)])]\, dr -\frac{1}{2} f_0\, [\bfcurY(r)]\, dr + \sqrt{f_0}\, [d\bfcurW(r)] \, ,                            
\end{equation}
with the initial condition
\begin{equation}
\label{EQ37}
[\bfcurZ(0)] = [\bfH_d]\, [a] \quad , \quad [\bfcurY(0)] = [\bfcurN\,] \, [a]\quad a.s \, ,                                                  
\end{equation}
in which the random matrices $[\curL([\bfcurZ(r)])]$ and $[d\bfcurW(r)]$ with values in $\MM_{\nu,m}$ are such that
\begin{equation}
\label{EQ38}
[\curL([\bfcurZ(r)])]= [L ( [\bfcurZ(r)] \, [g]^T ) ] \, [a]\, ,                                                                             
\end{equation}
\begin{equation}
\label{EQ39}
[d\bfcurW(r)]= [d\bfW(r)] \, [a]\, .                                                                                                  
\end{equation}
From Section~\ref{Section4.3}, it can be deduced that the problem defined by Eqs.~\eqref{EQ35} to  \eqref{EQ39} admits a unique invariant measure and a unique solution $\{ ([\bfcurZ(r)],[\bfcurY(r)]),$ $r\in \RR^+ \}$ that is a second-order diffusion stochastic process, which is stationary (for the shift semi-group on $\RR^+$) and ergodic.

For $\theta$ fixed in $\Theta$, the deterministic quantities $\{[\bfcurW(r;\theta)],r\geq 0\}$, $[\bfcurZ(0;\theta)]  = [\eta_d] \, [a]$, and $[\bfcurY(0;\theta)] = [\bfcurN(\theta)]\, [a]$  are independent realizations of the stochastic process $\{[\bfcurW(r)],r\geq 0\}$, the random matrix $[\bfcurZ(0)]$, and the random matrix $[\bfcurY(0)$.
Let $\{ ([\bfcurZ(r;\theta)],[\bfcurY(r;\theta)]), r\in \RR^+\}$ be the corresponding realization of the unique stationary diffusion process $\{ ([\bfcurZ(r)],[\bfcurY(r)]), r\in \RR^+\}$ of the reduced-order ISDE problem defined by Eqs.~\eqref{EQ35} to  \eqref{EQ37}). Then, using Eq.~\eqref{EQ28}, some additional realizations $[\eta_s^1], \ldots , [\eta_s^{n_\pMC}]$ of random matrix $[\bfH]$ can be generated by
\begin{equation}
\label{EQ40}
[\eta_s^\ell]  = [\bfcurZ(\ell \rho ;\theta)] \, [g]^T\quad , \quad \rho = M_0\, \Delta r \quad , \quad \ell = 1,\ldots n_\pMC\, ,               
\end{equation}
and using Eq.~\eqref{EQ3}, some additional realizations $[x_s^1], \ldots , [x_s^{n_\pMC}]$ of random matrix $[\bfX]$ can be generated (using the reduced-order representation defined by Eq.~\eqref{EQ3}) by
\begin{equation}
\label{EQ41}
[x_s^\ell] = [\underline x] + [\varphi]\, [\mu]^{1/2}\, [\eta_s^\ell]\quad , \quad \ell = 1,\ldots n_\pMC\, .                                  
\end{equation}
\subsection{Solving the reduced-order ISDE and computing the additional realizations for random matrix $[\bfX]$}
\label{Section4.7}
For numerically solving the reduced-order ISDE defined by Eqs.~\eqref{EQ35} to \eqref{EQ37}, a discretization scheme must be used.
For general surveys on discretization schemes for It\^o stochastic differential equations, we refer the reader to  \cite{Kloeden1992,Talay1990,Talay1995}. Concerning the particular cases related to Hamiltonian dynamical systems (which have also been analyzed in \cite{Talay2002} using an implicit
Euler scheme), we propose to use the St\"{o}rmer-Verlet scheme, which is a very efficient scheme that preserves energy for nondissipative Hamiltonian dynamical systems (see \cite{Hairer2002} for reviews about this scheme in the deterministic case, and see \cite{Burrage2007} and the references therein for the stochastic case).
\subsubsection{Discretization scheme of the reduced-order ISDE}
\label{Section4.7.1}
We then propose to reuse hereinafter the St\"{o}rmer-Verlet scheme, introduced and validated in \cite{Soize2012,Guilleminot2013,Soize2015}
for weakly dissipative stochastic Hamiltonian dynamical system.

Let $M = n_\MC \times M_0$ be the positive integer in which $n_\MC$ and $M_0$ have been introduced in Remark~5 of Section~\ref{Section4.3}. The reduced-order It\^o stochastic differential equation defined by Eqs.~\eqref{EQ35} and \eqref{EQ36} with the initial condition defined by Eq.~\eqref{EQ37}, is solved on the finite interval $\curR = [0\, , M\, \Delta r]$, in which $\Delta r$ is the sampling step of the continuous index parameter $r$. The integration scheme is based on the use of the $M+1$ sampling points $r_\ell$ such that $r_\ell = \ell\, \Delta r$ for $\ell =0,\ldots , M$.                                                             The following notations are introduced:
$[\bfcurZ_\ell] = [\bfcurZ(r_\ell)]$, $[\bfcurY_\ell] = [\bfcurY(r_\ell)]$, and $[\bfcurW_\ell] = [\bfcurW(r_\ell)]$, for $\ell =0,\ldots , M $,
 with
\begin{equation}
\label{EQ42}
[\bfcurZ_0] = [\bfH_d]\, [a]  \quad , \quad [\bfcurY_0] = [\bfcurN\,] \, [a] \quad , \quad [\bfcurW_0] = [0_{\nu,m}]  \quad a.s    \, .       
\end{equation}
For $\ell=0,\ldots , M-1$, let $[\Delta\bfcurW_{\ell+1}] = [\Delta\bfW_{\ell+1}]\, [a]$ be the sequence of random matrices with values in $\MM_{\nu,m}$, in which $[\Delta\bfW_{\ell+1}]= [\bfW_{\ell+1}] - [\bfW_\ell]$. The increments $[\Delta\bfW_{1}], \ldots , [\Delta\bfW_{M}]$ are $M$ independent random matrices.
For all $k=1,\ldots , \nu$ and for all $j=1,\ldots , N$, the real-valued random variables $\{[\Delta\bfW_{\ell+1}]_{kj}\}_{kj}$ are independent, Gaussian, second-order, and centered random variables such that
$E\{[\Delta\bfW_{\ell+1}]_{kj}[\Delta\bfW_{\ell+1}]_{k'j'}\}=\Delta r \,\delta_{kk'}\,\delta_{jj'}$.
For $\ell=0,\ldots, M-1$, the St\"{o}rmer-Verlet scheme applied to Eqs.~\eqref{EQ35} and \eqref{EQ36} yields
\begin{equation}
\label{EQ43}
[\bfcurZ_{\ell+\frac{1}{2}}]     =    [\bfcurZ_\ell] +\frac{\Delta r}{2} \, [\bfcurY_\ell] \, ,                                            
\end{equation}
\begin{equation}
\label{EQ44}
[\bfcurY_{\ell+1}]   =    \frac{1-b}{1+b}\, [\bfcurY_\ell] + \frac{\Delta r}{1+b}\, [\bfcurL_{\ell+\frac{1}{2}}] +
       \frac{\sqrt{f_0}}{1+b}\, [\Delta\bfcurW_{\ell+1}]\, ,                                                                              
\end{equation}
\begin{equation}
\label{EQ45}
[\bfcurZ_{\ell+ 1}]  =    [\bfcurZ_{\ell+\frac{1}{2}}] +\frac{\Delta r}{2} \, [\bfcurY_{\ell+1}]\, ,                                     
\end{equation}
with the initial condition defined by \eqref{EQ42}, where $b=f_0\, \Delta r\, / 4$, and where $[\bfcurL_{\ell+\frac{1}{2}}]$ is the $\MM_{\nu,m}$-valued random variable such that
\begin{equation}
\label{EQ46}
 [\bfcurL_{ \ell+\frac{1}{2} }]  = [\curL( [\bfcurZ_{\ell+\frac{1}{2}}] )] = [L( [\bfcurZ_{\ell+\frac{1}{2}}]\, [g]^T )]\, [a]\, ,        
 \end{equation}
in which, for all $[u] = [\bfu^1 \ldots \bfu^N]$ in $\MM_{\nu,N}$ with $\bfu^\ell=(u^\ell_1,\ldots ,u^\ell_\nu)$ in $\RR^\nu$, the entries of matrix $[L([u])]$ in $\MM_{\nu,N}$ are defined by Eqs.~\eqref{EQ17bis} and \eqref{EQ17ter}.
\subsubsection{Remarks about the estimation of the numerical integration parameters of the reduced-order ISDE}
\label{Section4.7.2}
Some estimations of the values of the parameters $f_0,\Delta r$, and $M_0$, which are used in the discretization scheme of the ISDE (with and without reduced-order representation of random matrix $[\bfH]$, and introduced in Sections~\ref{Section4.3} and \ref{Section4.7.1}) are described below.\\

\noindent (i) Parameter $\Delta r$ is written as $\Delta r = 2\pi\, \widehat s_\nu /\Fac$ in which $\Fac  > 1$ is an oversampling that has to be estimated  for getting a sufficient accuracy of the St\"{o}rmer-Verlet scheme (for instance, $\Fac =20$). This means that a convergence analysis of the solution must be carried out with respect to $\Fac$.\\

\noindent (ii) As the accuracy of the St\"{o}rmer-Verlet scheme  is finite, a small numerical integration error is unavoidably introduced. Although that the initial conditions are chosen in order to directly construct the stationary solution (associated with the unique invariant measure), a small transient response can occur and be superimposed to the stationary stochastic solution. Therefore, $f_0$ is chosen in order that  the damping in the dissipative Hamiltonian system is sufficiently large to rapidly kill such a small transient response (a typical value that is retained in the applications presented in Section~\ref{Section5} is $f_0=1.5$).\\

\noindent (iii) Using an estimation of the relaxation time of the underlying linear second-order dynamical system, and choosing an attenuation of $1/100$ for the transient response, parameter  $M_0$ must be chosen larger than
$2\log(100)\Fac/(\pi f_0 \, \widehat s_\nu)$.
A  typical value that is retained in the applications presented in Section~\ref{Section5} is $M_0 = 110$ or $330$).
\section{Applications}
\label{Section5}
Three applications are presented for random vector $\bfX$ with values in $\RR^n$ for which:
\begin{itemize}
\item  the dimension is $n=2$ and there are  $N=230$ given data points in subset $\curSn$, for which the mean value is made up of two circles in the plane).
\item the dimension is $n=3$ and there are $N=400$ given data points in subset $\curSn$, for which the mean value is made up of a helix in three-dimensional space).
\item the third example corresponds to a petro-physics database that
  is made up of experimental measurements (downloaded from
  \cite{boem}) and detailed in \cite{Thimmisetty2015}, for which the
  dimension is $n=35$ and for which $N=13,056$ given data points are
  concentrated in an unknown  "complex" subset $\curSn$ of $\RR^n$,
  which cannot be easily described once it is discovered.
\end{itemize}

\subsection{Application 1: Dimension $n=2$ with $N=230$ given data points}
\label{Section5.1}
\begin{figure}[!h]
\centering
\includegraphics[width=5.5cm]{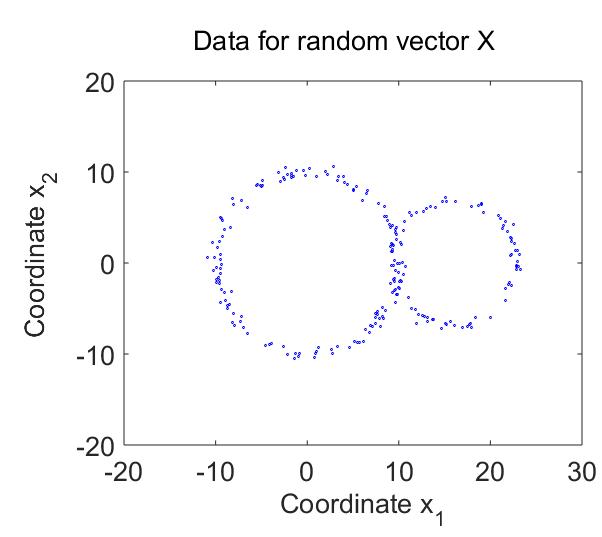} \hfil \includegraphics[width=5.5cm]{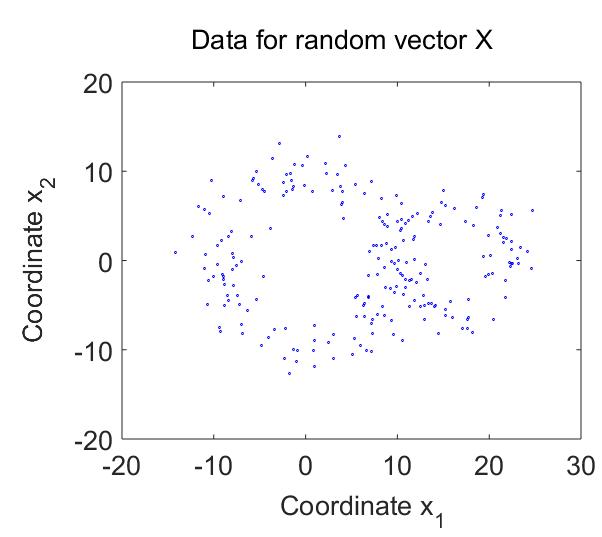}
\caption{$230$ given data points: case 1.1 (left), case 1.2 (right).}
\label{figure1}
\end{figure}

\begin{figure}[!h]
\centering
\includegraphics[width=5.5cm]{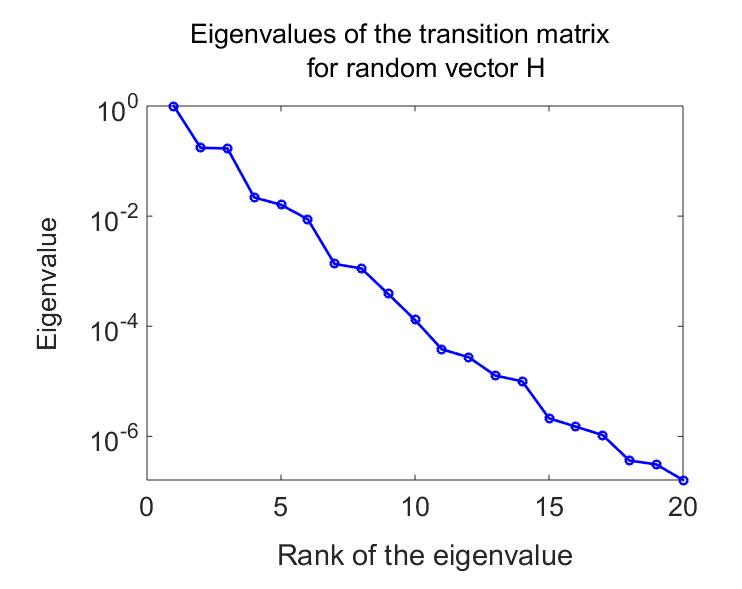} \hfil \includegraphics[width=5.5cm]{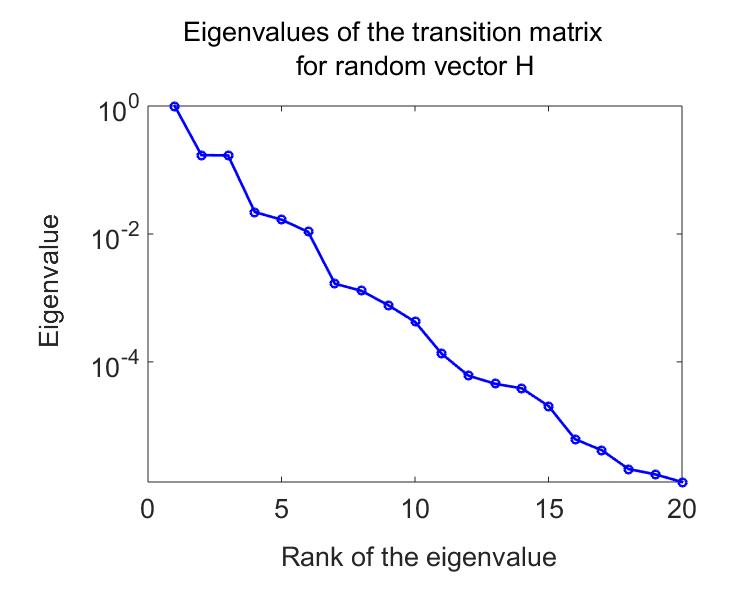}
\caption{Eigenvalues in $\log_{10}$-scale of the transition matrix for random vector $\bfH$: case 1.1 (left), case 1.2 (right).}
\label{figure2}
\end{figure}

\begin{figure}[!h]
\centering
\includegraphics[width=5.5cm]{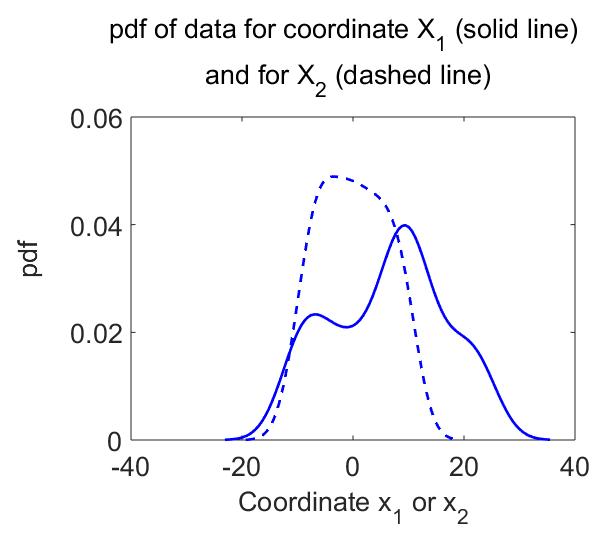} \hfil \includegraphics[width=5.5cm]{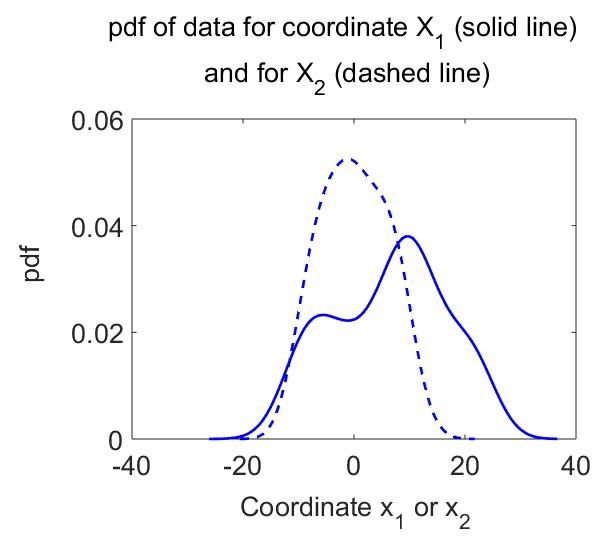}
\caption{pdf for random variables $X_1$ (solid line) and $X_2$ (dashed line) obtained by a nonparametric estimation from data points:
case 1.1 (left), case 1.2 (right).}
\label{figure3}
\end{figure}

\begin{figure}[!h]
\centering
\includegraphics[width=5.5cm]{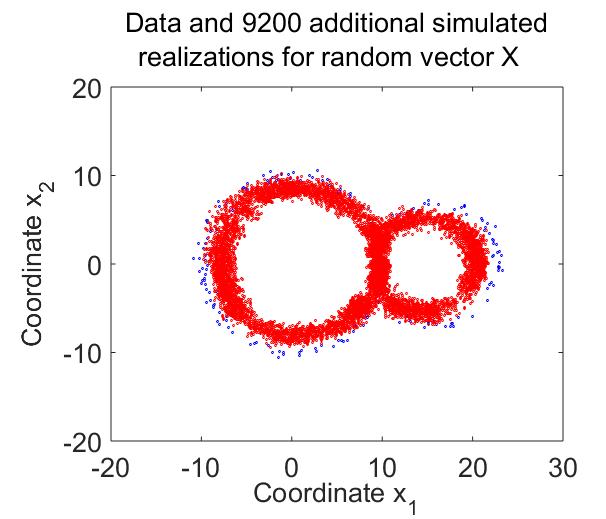} \hfil \includegraphics[width=5.5cm]{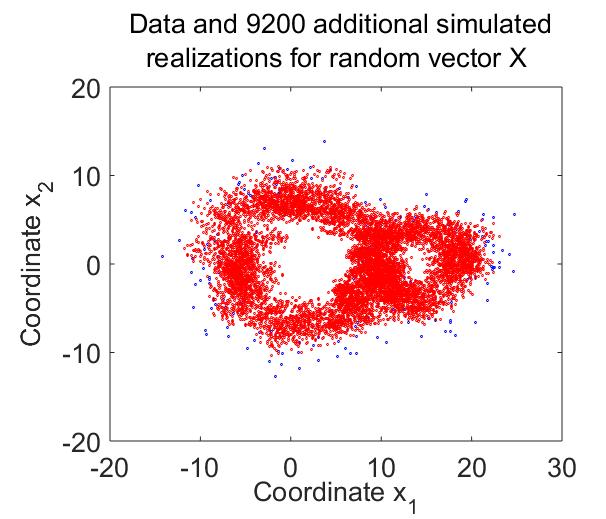}
\caption{$230$ given data points (blue symbols) and $9,200$ additional realizations (red symbols) generated using the reduced-order ISDE with $m=3$: case 1.1 (left), case 1.2 (right).}
\label{figure4}
\end{figure}

\begin{figure}[!h]
\centering
\includegraphics[width=5.5cm]{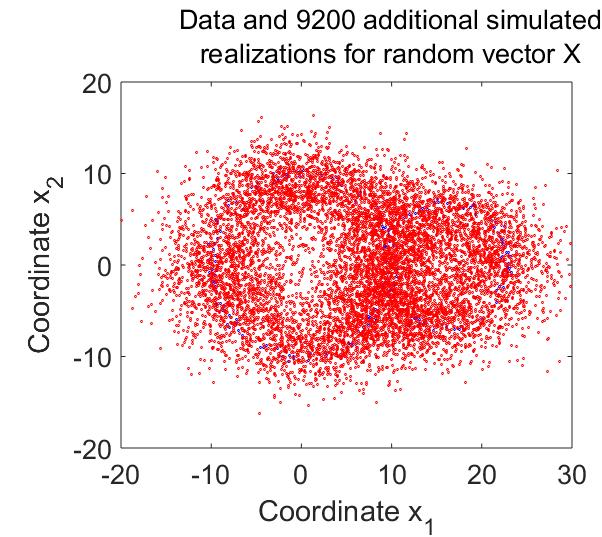} \hfil \includegraphics[width=5.5cm]{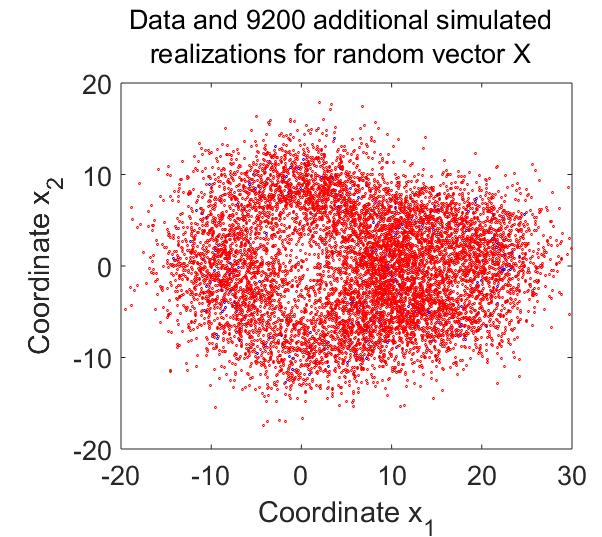}
\caption{$230$ given data points (blue symbols) and $9,200$ additional realizations (red symbols) generated using the ISDE: case 1.1 (left), case 1.2 (right).}
\label{figure5}
\end{figure}

For this first application, two cases are considered: small (case 1.1) and medium (case 1.2) statistical fluctuations around the two circles. For every case, the number of given data points is $N=230$, no scaling of data is performed, but the normalization defined in Section~\ref{Section4.1} is done and yields $\nu=2$.
In Figs.~\ref{figure1} to \ref{figure5}, the left figures are relative to case 1.1 and the right ones to case 1.2.
Fig.~\ref{figure1} displays the $230$ given data points for random vector $\bfX =(X_1,X_2)$ of the data set represented by matrix $[x^d]$ in $\MM_{2,230}$, and shows that the given data points are concentrated in the neighborhood of two circles, with small  (case 1.1) and medium (case 1.2) statistical fluctuations.
The kernel is defined by Eq.~\eqref{EQ19}, the value of the smoothing parameter that is retained is $\varepsilon=2.7318$, $\kappa$ is chosen to $1$, and the graph of the eigenvalues of the transition matrix for random vector $\bfH$ is displayed in Fig.~\ref{figure2}. These two graphs show that dimension $m$ can be chosen to $3$, and for $m=3$, the value of $e_\red(m)$ (defined by Eq.~\eqref{EQ34}) is $6.34\times 10^{-4}$ for case 1.1 and $9.28\times 10^{-4}$ for case 1.2. It can thus be considered that a reasonable mean-square convergence is reached for these two cases.
Fig.~\ref{figure3} displays the pdf for random variables $X_1$ and $X_2$  computed with a nonparametric estimation from the data points.
For all the computation, the numerical values of the parameters for generating $9,200$ additional realizations are $\Delta r=0.1179$, $M_0=110$, and $n_\MC = 40$, yielding $M=4,400$.
The results obtained with the reduced-order ISDE (for which the first $m=3$ vectors of the diffusion-maps basis are used) are displayed in
Fig.~\ref{figure4}, which shows the $230$ given data points and the $9,200$ additional realizations generated using the reduced-order ISDE. It can be seen that the additional realizations are effectively concentrated in subset $\curSn$.
Fig.~\ref{figure5} displays the $230$ given data points  and the $9,200$ additional realizations  generated using a direct simulation of the ISDE  presented in Section~\ref{Section4.3}.
It can be seen that the realizations are not concentrated in subset $\curSn$, but are scattered.

\subsection{Application 2: Dimension $n=3$ with $N=400$ given data points}
\label{Section5.2}
\begin{figure}[!h]
\centering
\includegraphics[width=5.5cm]{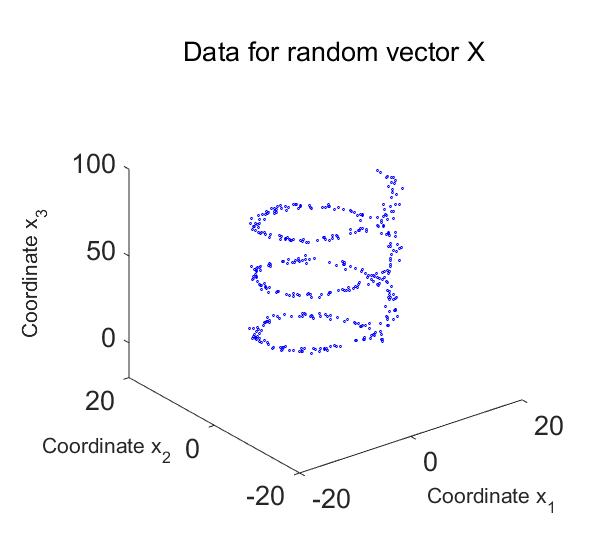} \hfil \includegraphics[width=5.5cm]{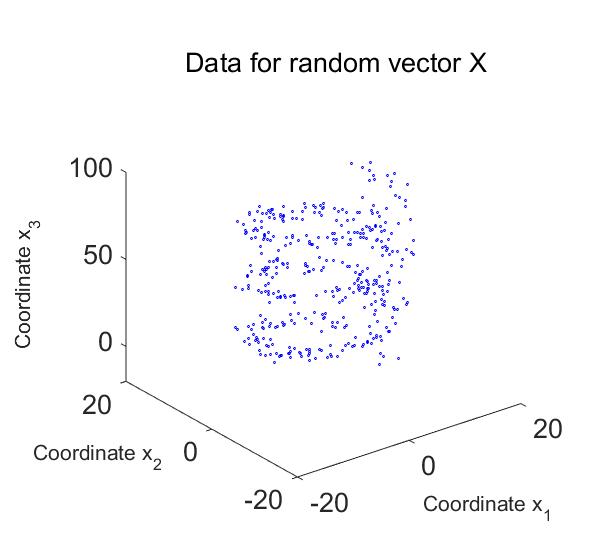}
\caption{$400$ given data points: case 2.1 (left), case 2.2 (right).}
\label{figure6}
\end{figure}

\begin{figure}[!h]
\centering
\includegraphics[width=5.5cm]{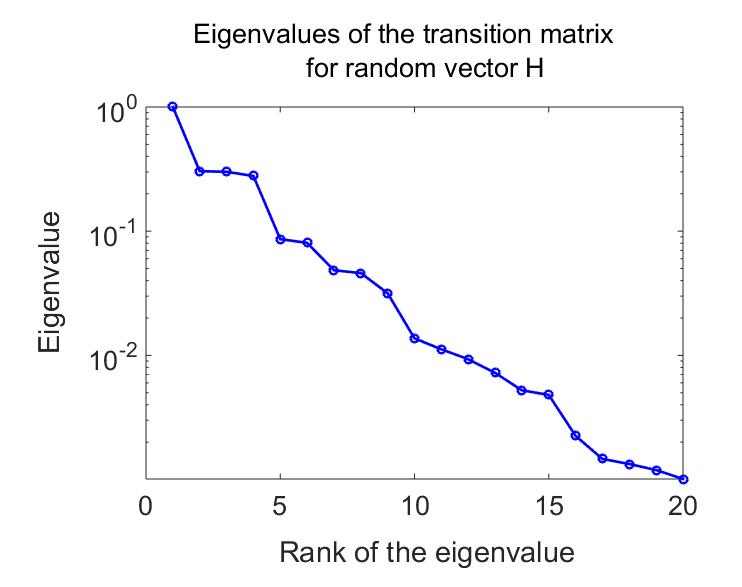} \hfil \includegraphics[width=5.5cm]{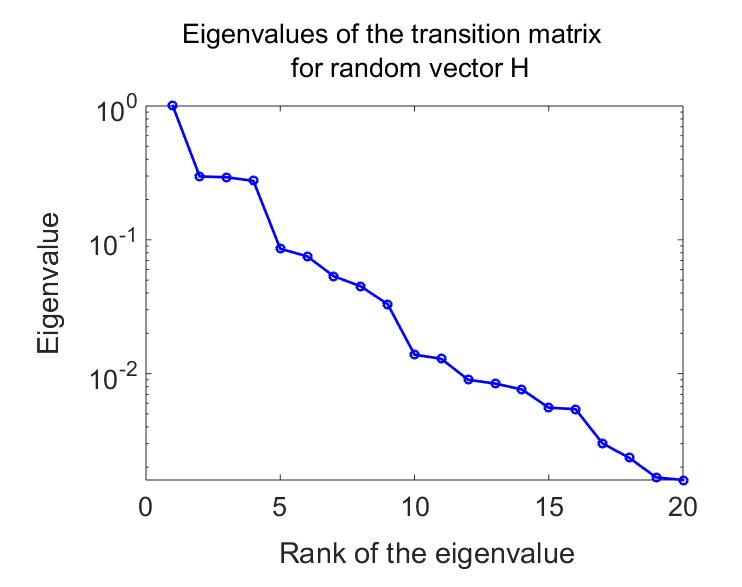}
\caption{Eigenvalues in $\log_{10}$-scale of the transition matrix for random vector $\bfH$: case 2.1 (left), case 2.2 (right).}
\label{figure7}
\end{figure}

\begin{figure}[!h]
\centering
\includegraphics[width=5.5cm]{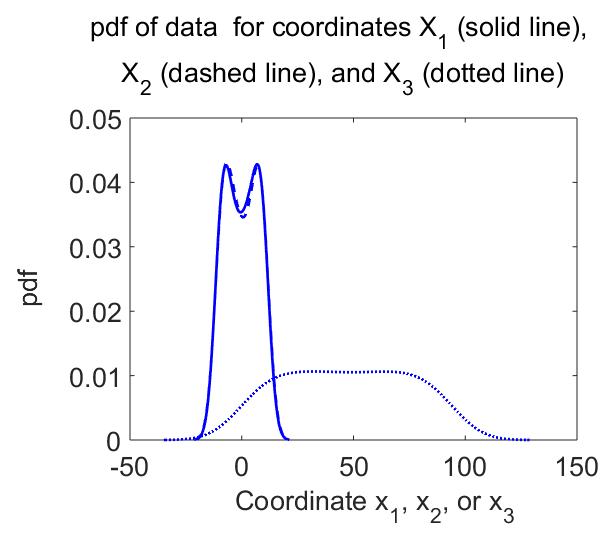} \hfil \includegraphics[width=5.5cm]{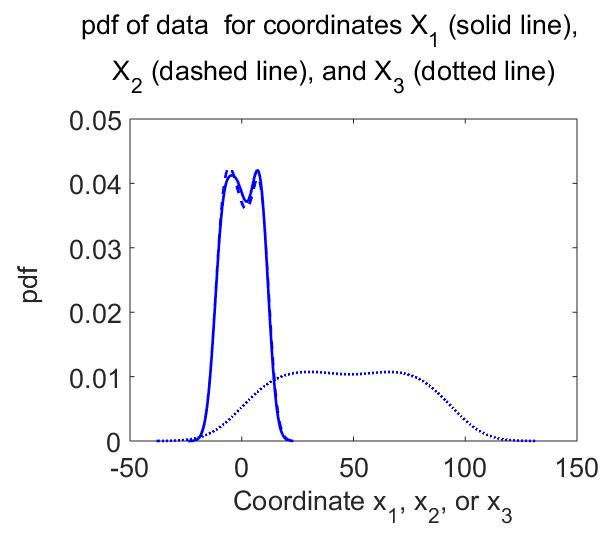}
\caption{pdf for random variables $X_1$ (solid line), $X_2$ (dashed line), and $X_3$ (dotted line) obtained by a nonparametric estimation from data points:
case 2.1 (left), case 2.2 (right).}
\label{figure8}
\end{figure}

\begin{figure}[!h]
\centering
\includegraphics[width=5.5cm]{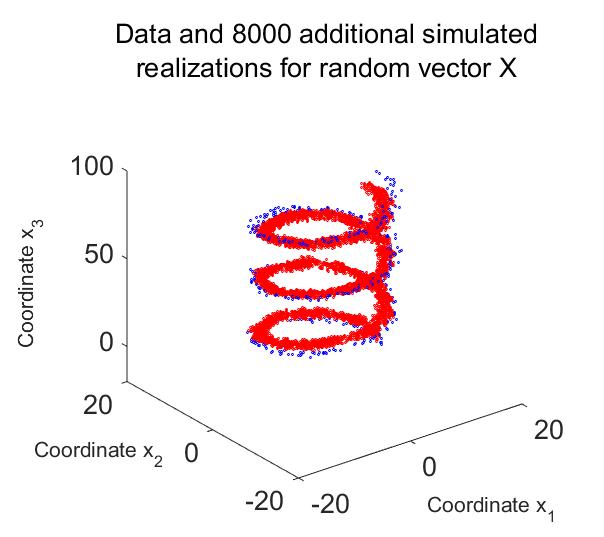} \hfil \includegraphics[width=5.5cm]{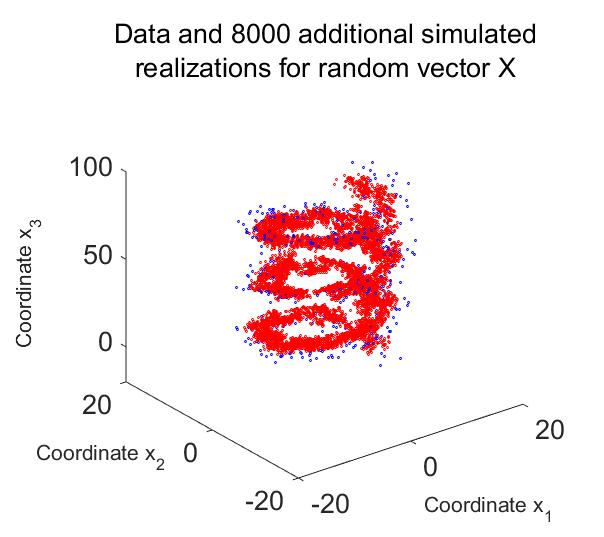}
\caption{$400$ given data points (blue symbols) and $8,000$ additional realizations (red symbols) generated using the reduced-order ISDE with $m=4$: case 2.1 (left), case 2.2 (right).}
\label{figure9}
\end{figure}

\begin{figure}[!h]
\centering
\includegraphics[width=5.5cm]{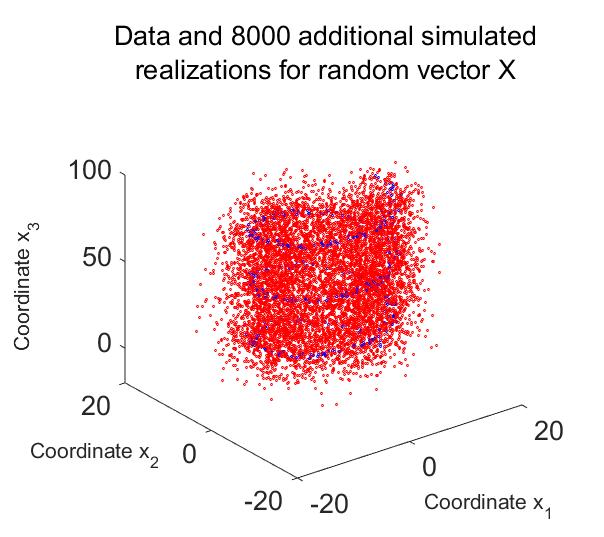} \hfil \includegraphics[width=5.5cm]{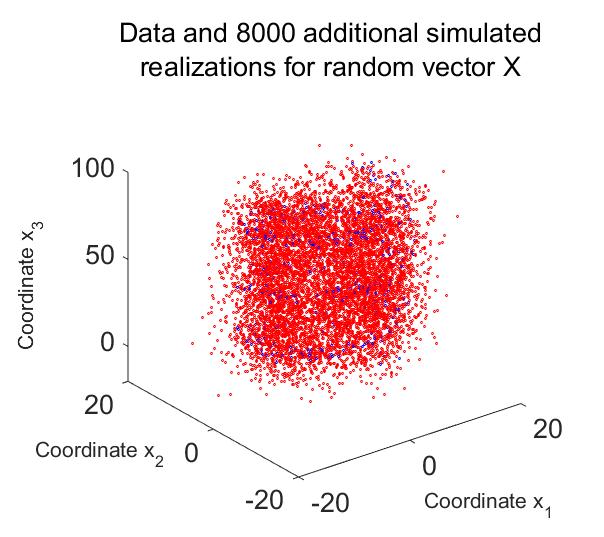}
\caption{$400$ given data points (blue symbols) and $8,000$ additional realizations (red symbols) generated using the ISDE: case 2.1 (left), case 2.2 (right).}
\label{figure10}
\end{figure}

As previously, two cases are considered: small (case 2.1) and medium (case 2.2) statistical fluctuations around the helical. For every case, the number of given data points is $N=400$, no scaling of data is performed, but the normalization defined in Section~\ref{Section4.1} is done and yields $\nu=3$.
In Figs.~\ref{figure6} to \ref{figure10}, the left figures are relative to case 2.1 and the right ones to case 2.2.
Fig.~\ref{figure6} displays the $400$ given data points for random vector $\bfX =(X_1,X_2,X_3)$ of the data set represented by matrix $[x^d]$ in $\MM_{3,400}$. Fig.~\ref{figure6} shows that the given data points are concentrated in the neighborhood of the helical, with small  (case 2.1) and medium (case 2.2) statistical fluctuations.
The kernel is defined by Eq.~\eqref{EQ19}, the value of the smoothing parameter that is retained is $\varepsilon=1.57$, $\kappa$ is chosen to $1$, and the graph of the eigenvalues of the transition matrix for random vector $\bfH$ is displayed in Fig.~\ref{figure7}. These two graphs show that dimension $m$ can be chosen to $4$, and for $m=4$, the value of $e_\red(m)$ (defined by Eq.~\eqref{EQ34}) is $5.53\times 10^{-4}$ for case 2.1 and $4.28\times 10^{-4}$ for case 2.2. It can thus be considered that a reasonable mean-square convergence is reached for these two cases.
Fig.~\ref{figure8} displays the pdf for random variables $X_1$, $X_2$, and $X_3$  computed with a nonparametric estimation from the data points.
For all the computation, the numerical values of the parameters for generating $9,200$ additional realizations are $\Delta r=0.1196$, $M_0=110$, and $n_\MC = 20$, yielding $M=2,200$.
The results obtained with the reduced-order ISDE (for which the first $m=4$ vectors of the diffusion-maps basis are used) are displayed in
Fig.~\ref{figure9}, which shows the $400$ given data points and the $8,000$ additional realizations generated using the reduced-order ISDE. It can be seen that the additional realizations are effectively concentrated in subset $\curSn$.
Fig.~\ref{figure10} displays the $400$ given data points  and the $8,000$ additional realizations  generated using a direct simulation with the ISDE presented in Section~\ref{Section4.3}.
It can be seen that the realizations are not concentrated in subset $\curSn$, but are scattered.
\subsection{Application 3: Dimension $n=35$ with $N=13,056$ given data points}
\label{Section5.3}
\begin{figure}[!ht]
\centering
\includegraphics[width=4.4cm]{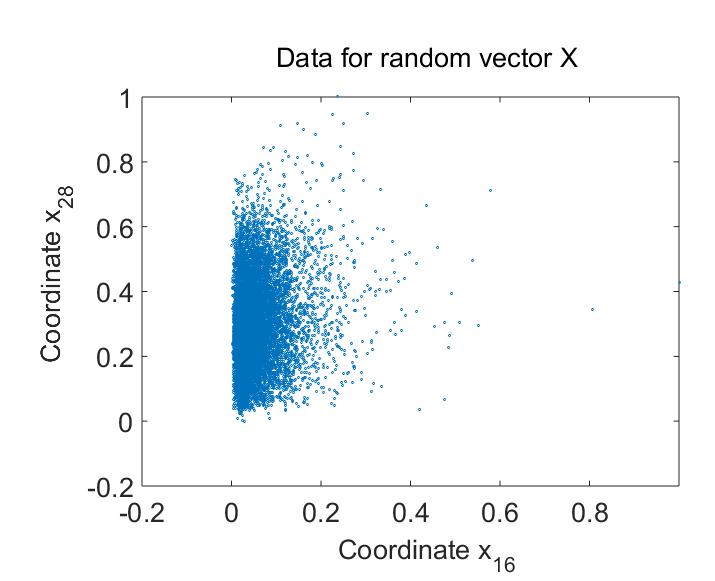}  \includegraphics[width=4.2cm]{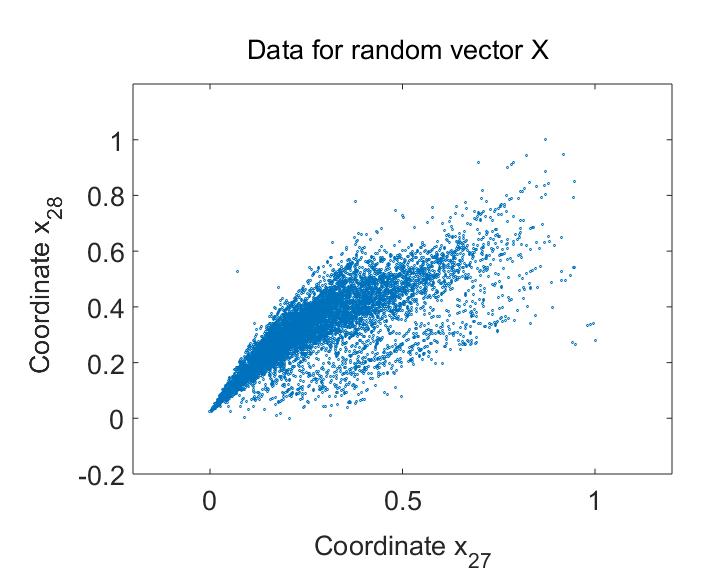}  \includegraphics[width=4.8cm]{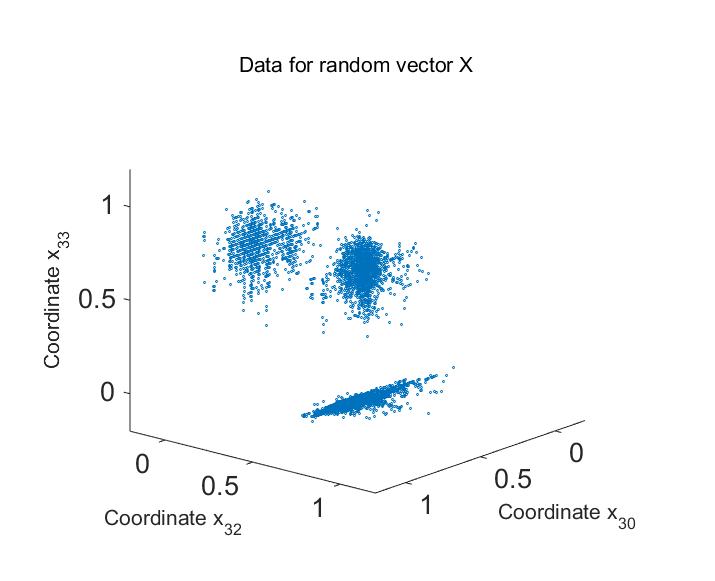}
\caption{$13,056$ given data points viewed from coordinates $x_{16}$ and $x_{28}$ (up left), viewed from coordinates $x_{27}$ and $x_{28}$ (up right),
and viewed from coordinates $x_{30}$, $x_{32}$, and $x_{33}$ (down).}
\label{figure11}
\end{figure}
\begin{figure}[!ht]
\centering
\includegraphics[width=5.0cm]{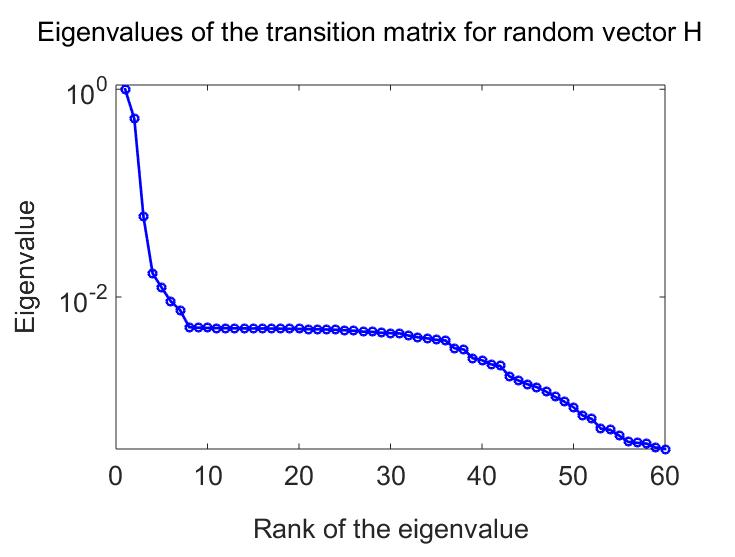} \hfil \includegraphics[width=5.0cm]{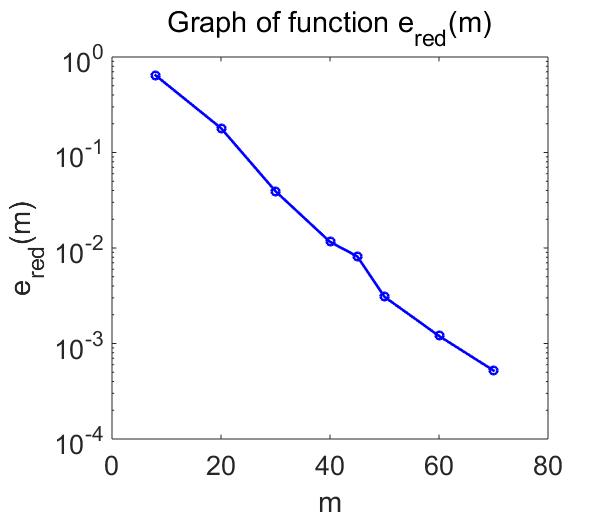}
\caption{Eigenvalues in $\log_{10}$-scale of the transition matrix for random vector $\bfH$ (left). Graph $m\mapsto e_{\red}(m)$ in $\log_{10}$ scale (right).}
\label{figure12}
\end{figure}
\begin{figure}[!ht]
\centering
\includegraphics[width=5.5cm]{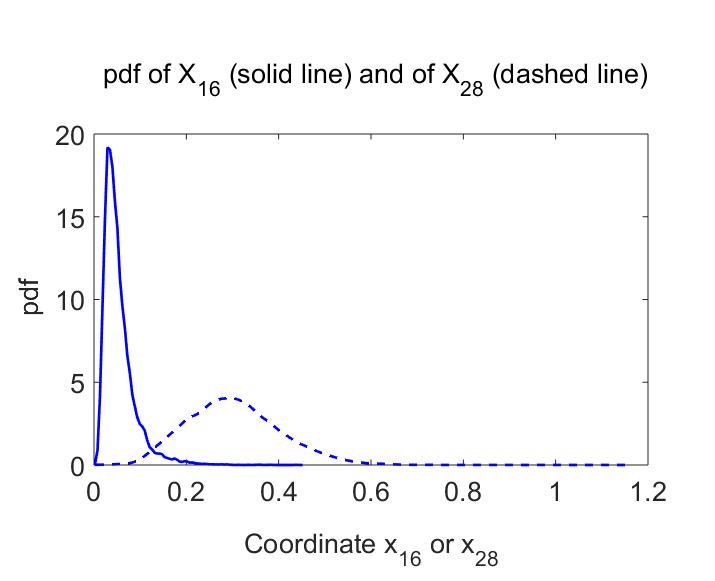} \hfil \includegraphics[width=5.0cm]{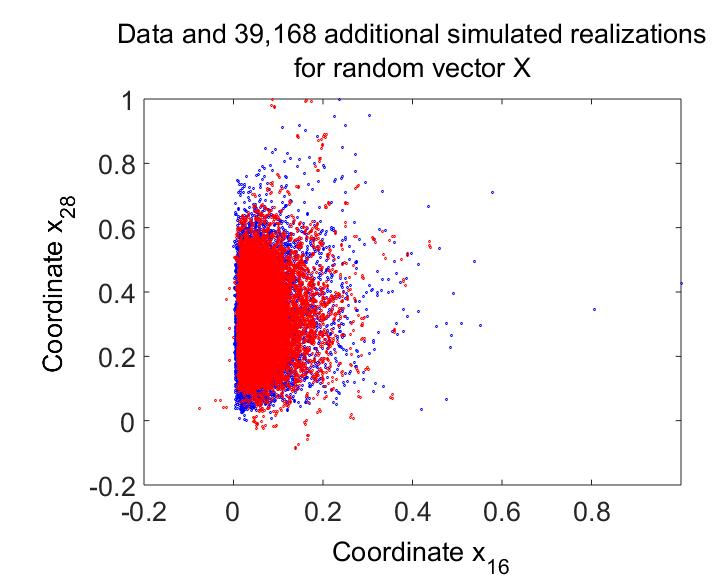}
\\

\includegraphics[width=5.5cm]{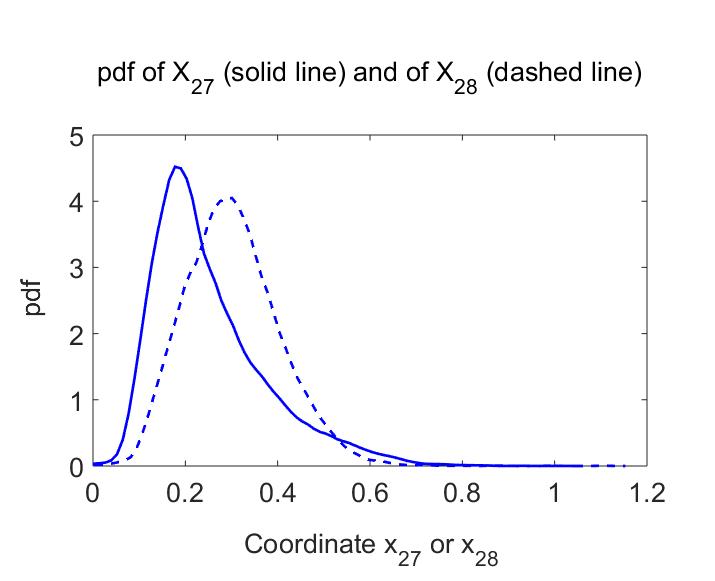} \hfil \includegraphics[width=5.0cm]{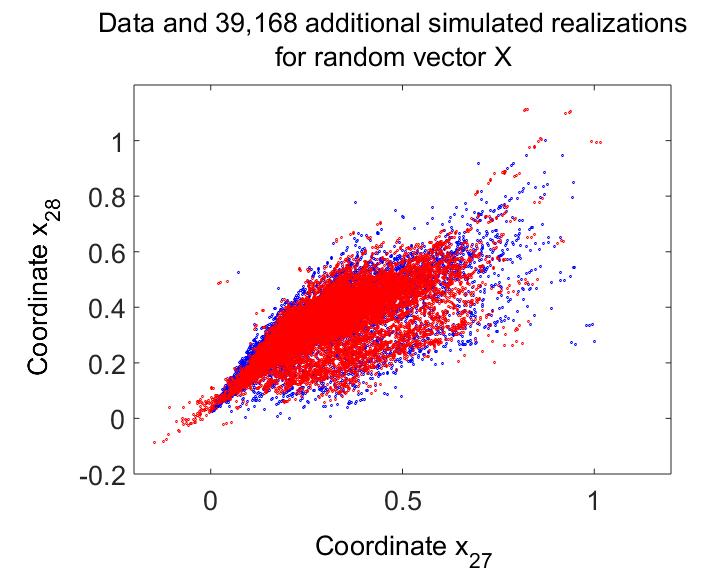}
\\

\includegraphics[width=5.5cm]{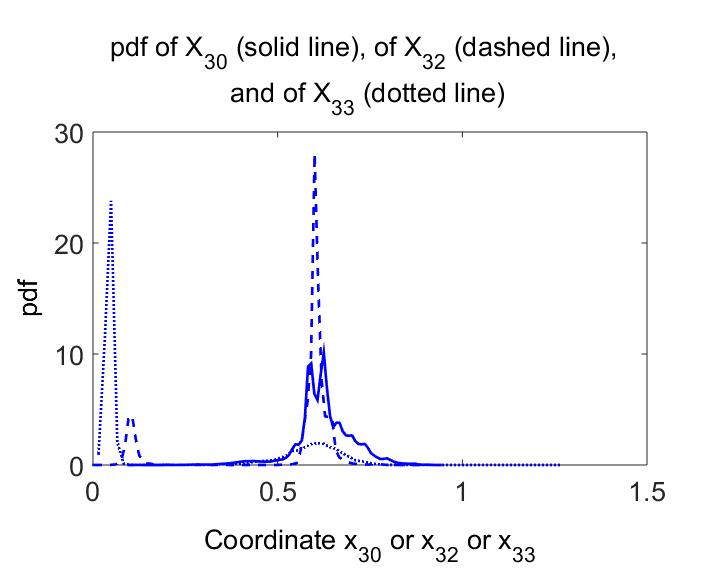} \hfil \includegraphics[width=5.0cm]{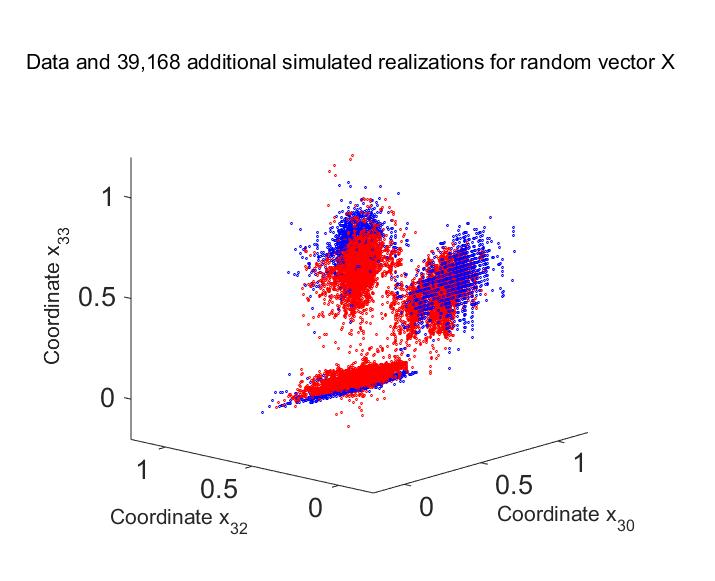}

\caption{Left figures: Illustration of the pdf for some components of random vector $\bfX$ obtained by a nonparametric estimation from the data points and the simulated data points. Right figures: $13,056$ given data points (blue symbols) and $39,168$ additional realizations (red symbols) generated using the reduced-order representation of $[\bfH]$  with $m=50$,  viewed from different components of random vector $\bfX$.}
\label{figure13}
\end{figure}
\begin{figure}[!ht]
\centering
\includegraphics[width=4.4cm]{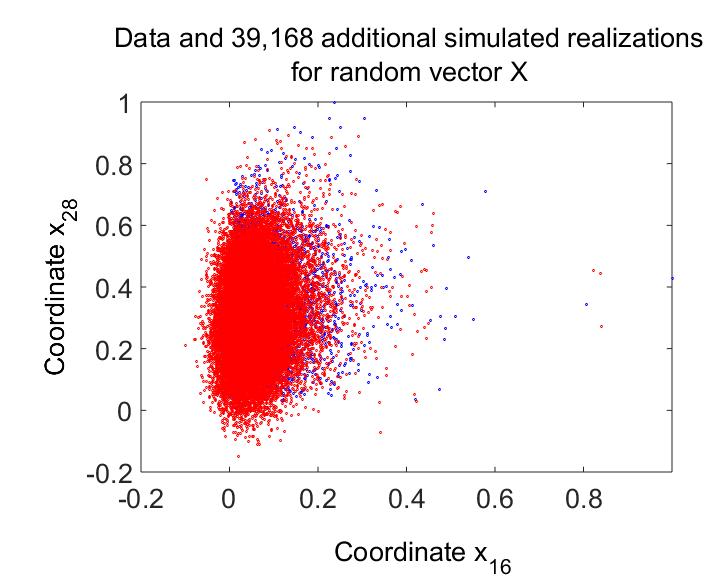}  \includegraphics[width=4.2cm]{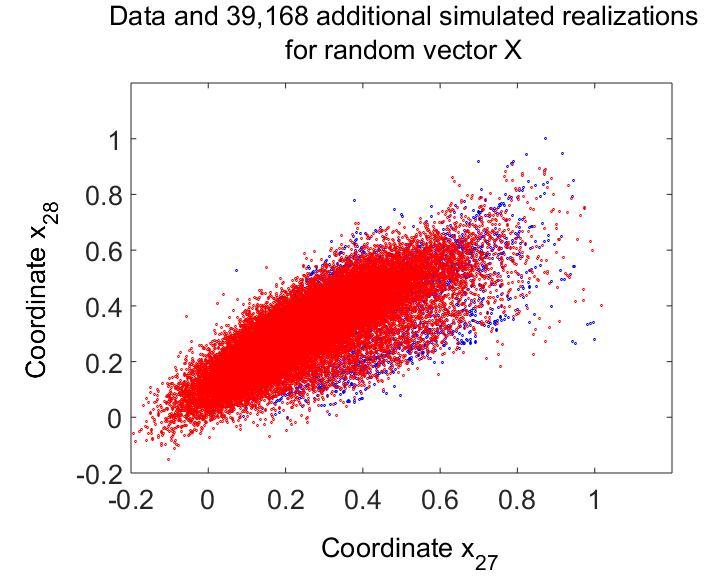} \includegraphics[width=4.8cm]{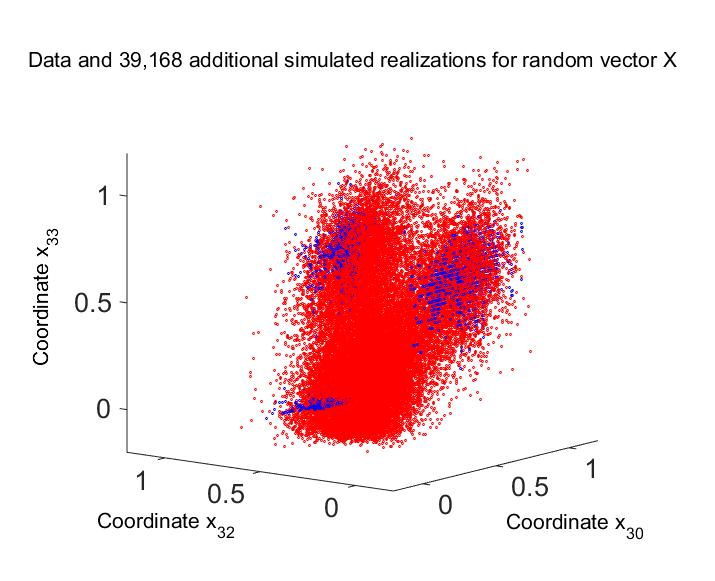}
\caption{$13,056$ given data points (blue symbols) and $39,168$ additional realizations (red symbols) generated  without using the reduced-order representation,  viewed from different components of random vector $\bfX$.}
\label{figure15}
\end{figure}
The data base used corresponds to a petro-physics data base of experimental experiments. The dimension of random vector $\bfX$ is $n=35$ and the number of given data points is $N=13,056$. The scaling and the normalization defined in Section~\ref{Section4.1} are necessary, have been done, and yield $\nu=32$.
Fig.~\ref{figure11} displays $13,056$ given data points viewed from coordinates $x_{16}$ and $x_{28}$, from coordinates $x_{27}$ and $x_{28}$,
and from coordinates $x_{30}$, $x_{32}$, and $x_{33}$. Although only a partial representation of the $13,056$ data points for the $\RR^n$-valued random vector $\bfX$ is given, this figure shows that $\curSn$ is certainly a complex subset of $\RR^n$.
The kernel is defined by Eq.~\eqref{EQ19}, the value of the smoothing parameter that has been used is $\varepsilon=100$, and $\kappa$ has also been chosen to $1$. The graph of the eigenvalues (of the transition matrix relative to random vector $\bfH$) displayed in Fig.~\ref{figure12} (left) shows that the value $m=8$ could potentially be a good choice for the value of $m$. However, for $m=8$, the value of $e_\red(m)$ is $0.99$ that shows that the mean-square convergence is not reached. Consequently, an analysis has been performed in constructing the graph of function $m\mapsto e_\red(m)$ in order to identify the smallest value of $m$ for which the mean-square convergence is reasonably reached. The graph displays in Fig.~\ref{figure12} (right) clearly shows that a good choice is $m=50$ for which the value of $e_\red(m)$ is $3.08\times 10^{-3}$ that can thus be considered as a reasonable mean-square convergence.
For all the computation, the numerical values of the parameters for generating $39,168$ additional realizations are $\Delta r= 0.06142$, $M_0=330$, and $n_\MC = 3$, yielding $M=990$.\\
For the same coordinates that those introduced in Fig.~\ref{figure11}, the left figures in Fig.~\ref{figure13} display the pdf of the considered components of random vector $\bfX$ obtained by a nonparametric estimation from the data points and the simulated data points obtained with the reduced-order ISDE, and the right figures display the $13,056$ given data points and the $39,168$ additional realizations generated using the reduced-order ISDE using the first $m=50$ vectors of the diffusion-maps basis. It can be seen that the additional realizations are effectively concentrated in subset $\curSn$.
Fig.~\ref{figure15} displays the $13,056$ given data points  and the $39,168$ additional realizations  generated using a direct simulation with the ISDE presented in Section~\ref{Section4.3}.
It can be seen that the realizations are not concentrated in subset $\curSn$, but are scattered. In particular, the positivity of random variable $X_{16}$ is not satisfied.

\section{Conclusions}
\label{Section6}
A new methodology has been presented and validated for generating
realizations of an $\RR^n$-valued random vector, for which the
probability distribution is unknown and is concentrated on an unknown
subset $\curSn$ of $\RR^n$. Both the probability distribution and the subset
$\curSn$ are constructed to be statistically consistent with a
specified data set construed as providing initial realizations of the
random vector.  The proposed method is robust and can be used for high
dimension and for large initial data sets. It is expected that  the proposed
method will contribute to open new possibilities of developments in
many areas of uncertainty quantification and statistical data
analysis, in particular in the design of experiments for random
parameters.

\section{Acknowledgment}

Part of this research was supported by the U.S. Department of Energy
Office of Advanced Scientific Computing Research.





%
%

\end{document}